\documentstyle{article}

\newcommand{\be}{\begin{equation}}
\newcommand{\ee}{\end{equation}}
\newcommand{\ba}{\begin{eqnarray}}
\newcommand{\ea}{\end{eqnarray}}
\newcommand{\baa}{\begin{eqnarray*}}
\newcommand{\eaa}{\end{eqnarray*}}
\newcommand{\bb}{\begin{thebibliography}}
\newcommand{\eb}{\end{thebibliography}}

\newcommand{\bi}[1]{\bibitem{#1}}
\newcommand{\lab}[1]{\label{#1}}
\newcommand{\re}[1]{(\ref{#1})}

\newcounter{my}
\newcommand{\he}%
   {\stepcounter{equation}\setcounter{my}%
   {\value{equation}}\setcounter{equation}0%
   }%
\newcommand{\she}%
   {\setcounter{equation}{\value{my}}%
    }%

\renewcommand\t{\tilde}

\newcommand\ve{\varepsilon}

\newtheorem{lem}{Lemma}

\begin{document}

%\begin{titlepage}

\begin{center}

{\Large \bf Dunkl shift operators and Bannai-Ito polynomials}

\vspace{5mm}

{\large \bf Satoshi Tsujimoto}

\medskip

{\it Department of Applied Mathematics and Physics, Graduate
School of Informatics, Kyoto University, Sakyo-ku, Kyoto
606--8501, Japan}

\medskip

\vspace{5mm}

{\large \bf Luc Vinet}

\medskip

{\it Centre de Recherches Math\'ematiques Universit\'e de
Montr\'eal, P.O. Box 6128, Centre-ville Station, Montr\'eal
(Qu\'ebec), H3C 3J7, Canada}

\medskip

and

\medskip

{\large \bf Alexei Zhedanov}

\medskip

{\em Donetsk Institute for Physics and Technology, Donetsk 83114,  Ukraine}

\end{center}

%\vspace*{5mm}

\begin{abstract}
We consider the most general Dunkl shift operator $L$ with the
following properties: (i) $L$ is of first order in the shift
operator and involves reflections; (ii) $L$ preserves the space of
polynomials of a given degree; (iii) $L$ is potentially
self-adjoint. We show that under these conditions, the operator
$L$ has eigenfunctions which coincide with the Bannai-Ito
polynomials. We construct a polynomial basis which is
lower-triangular and two-diagonal with respect to the action of
the operator $L$. This allows to express the BI polynomials
explicitly. We also present an anti-commutator AW(3) algebra
corresponding to this operator. From the representations of this
algebra, we derive the structure and recurrence relations of the
BI polynomials. We introduce new orthogonal polynomials - referred
to as the complementary BI polynomials - as an alternative $q \to
-1$ limit of the Askey-Wilson polynomials. These complementary BI
polynomials lead to a new explicit expression for the BI
polynomials in terms of the ordinary Wilson polynomials.

\vspace{2cm}


\vspace{2cm}


\end{abstract}

%\end{titlepage}

\newpage
\section{Introduction}
Two new families of ``classical`` polynomials \cite{VZ_little},
\cite{VZ_big} were introduced recently through limits when $q$
goes to -1 of the little and big q-Jacobi polynomials.  The term
"classical" is taken to mean that the polynomials $P_n(x)$ are
eigenfunctions of some differential or difference operator $L$:
\be L P_n(x) = \lambda_n P_n(x). \lab{LPP0} \ee In these
instances, the operator $L$ is of first order in the derivative
operator $\partial_x$ and contains moreover the reflection
operator $R$ defined by $R f(x)=f(-x)$; it can be identified as a
first order operator of Dunkl type \cite{Dunkl} written as \be L =
F(x)(I-R) + F_1(x)
\partial_x + G(x)\partial_x R  \lab{L_Dunkl} \ee
with some real rational functions $F(x), F_1(x), G(x)$.

The eigenvalue equation \re{LPP0} (with $L$ as in \re{L_Dunkl})
has further been investigated generally in \cite{VZ_Bochner}. As a
rule the polynomial eigenfunctions $P_n(x)$ are not orthogonal.
However, if $F_1(x)=0$ the big and little -1 Jacobi polynomials
are seen to be the only classes of orthogonal polynomials
consisting in eigenfunctions of such first order Dunkl
differential operators.

The main purpose of the present paper is to study a difference
analogue of the operator $L$ and the corresponding
eigenpolynomials $P_n(x)$. We introduce a natural generalization
of the operator $L$ which (apart from the multiplication operator)
only contains 2 nontrivial  operators: the shift $T^{+}$ and the
reflection $R$. We shall call such operators Dunkl shift operators
on the uniform grid. We find the conditions for such operators to
transform any polynomial into a polynomial of the same degree.
Under these conditions  the eigenvalue problem can be posited and
we can look for the corresponding eigenpolynomial solutions
$P_n(x)$.

The main result of our study is that the polynomial eigenfunctions
of such a generic Dunkl shift operator, coincide with the
Bannai-Ito (BI) polynomials that depend on 4 parameters.

Recall that the Bannai-Ito polynomials were first proposed in
\cite{BI} as a  $q \to -1$ limit of the q-Racah polynomials.
Bannai and Ito showed that these polynomials together with the
q-Racah polynomials (and their specializations and limiting cases)
are the most general orthogonal polynomial systems satisfying the
Leonard duality property. Bannai and Ito derived the three-term
recurrence relation for their polynomials and also presented
explicit expressions in terms of linear combinations of two
hypergeometric functions ${_4}F_3(1)$ \cite{BI}.

Further developments of the theory of BI polynomials and their
applications can be found in papers by Terwilliger \cite{Ter},
\cite{Ter2}, Curtin \cite{Curtin1}, \cite{Curtin2} and Vid\=unas
\cite{Vidunas}.

Our approach is completely different: we start with a generic
first order Dunkl shift operator that preserves the space of
polynomials and then show that its polynomial eigensolutions
coincide with the BI polynomials. This approach gives a new type
of eigenvalue equation for the BI polynomials  - one that involves
a combination of the ordinary shift operator with reflections.
Moreover, we give an algebraic interpretation of the BI
polynomials in terms of a $q=-1$ version of the AW(3) algebra.
This interpretation allows one to derive the structure and
recurrence relations for the BI polynomials. We introduce also a
new class of orthogonal polynomials - the complementary Bannai-Ito
polynomials. They can be obtained from the Askey-Wilson
polynomials by a limiting process $q \to -1$ which is slightly
different from the one used by Bannai and Ito in \cite{BI} to
obtain their polynomials. In contrast to the BI polynomials, the
complementary BI polynomials do not possess the Leonard duality
property. Nevertheless, the complementary BI polynomials have a
very simple expression in terms of the ordinary Wilson
polynomials. Moreover, the complementary BI polynomials are
Christoffel transforms of the BI polynomials. This allows one to
present a new explicit expression of the BI polynomials in the
form of a linear combination of two Wilson polynomials.

The present paper is organized as follows.

In the next section we introduce a class of (potentially
symmetric) operators $L$ which has the following properties:

(i) $L$ is of first order with respect to $T^{+}$ and contains
also the reflection operator $R$;

(ii) $L$ preserves the space of polynomials (i.e. it transforms
any polynomial into a polynomial of the same degree).

These properties imply the existence of eigenpolynomials $P_n(x)$
satisfying \re{LPP0}.

In Section 3, we determine these eigenpolynomials $P_n(x)$
explicitly. Most useful to that end is the observation that there
is a polynomial basis ($\phi_n(x)$) in which the operator $L$ is
lower triangular with only 2 diagonals, i.e.
$$
L \phi_{n}(x)=\lambda_n \phi_n(x) + \nu_n \phi_{n-1}(x) .
$$
This property leads to an explicit formula for the polynomials
$P_n(x)$ in terms of a linear combination of two hypergeometric
functions ${_4}F_3(1)$.

In Section 4, we show that the operator $L$ together with the
operator multiplication by $x$, forms a $q=-1$ (i.e. an
anticommutator) analogue of the $AW(3)$ algebra which we call the
Bannai-Ito algebra. Using the representations of this algebra, we
determine the structure relations of the BI polynomials $P_n(x)$
and  show that these polynomials satisfy expectedly a 3-term
recurrence relation
$$
P_{n+1}(x) + b_n P_n(x) + u_n P_{n-1}(x) = xP_n(x),
$$
where the recurrence coefficients $u_n, b_n$ are derived
explicitly. This allows to identify the eigenpolynomials $P_n(x)$
with the Bannai-Ito polynomials first introduced in \cite{BI}.

In Section 5, we construct the complementary Bannai-Ito
polynomials and obtain a new explicit expression of the BI
polynomials in terms of ordinary Wilson polynomials.

In Section 6, we demonstrate that the operator $L$ is
symmetrizable. This means that there exists a function
$\varphi(x)$ such that the operator $\varphi(x) L$ is symmetric
$$
(\varphi(x) L)^* = \varphi(x) L,
$$
where $L^*$ means the formal conjugate of $L$ (constructed
according to natural rules).

In Section 7, we introduce the Bannai-Ito grid, which is a
discrete invariant set of the BI operator. The BI operator becomes
a 3-diagonal matrix when restricted to this set. Using this
observation we derive the weight function for the case when the BI
polynomials are orthogonal on a finite set of points of the real
line.

In Section 8, it is shown how the BI and complementary BI
polynomials can be obtained by limiting processes from the
Askey-Wilson polynomials.

In Section 9, it is indicated how the Dunkl shift operator can be
derived from the Askey-Wilson difference operator when $q \to -1$.

In Section 10, we consider a special symmetric case of the BI
polynomials. This leads to the construction of a new  family of
orthogonal polynomials (explicitly expressed in terms of dual
continuous Hahn polynomials) with a purely continuous weight
function on the whole real axis. This is the first nontrivial
example of an infinite positive family of BI polynomials. (Recall
that only finite systems of BI polynomials have been considered in
details so far).

In Section 11 we show that a limiting process $x \to x/h, h \to 0$
reduces the Dunkl shift operator to a Dunkl differential operator
$L$ having the form \re{L_Dunkl} with $F_1(x)=0$. This means that
correspondingly, the  Bannai-Ito polynomials tend in this limit to
the big -1 Jacobi polynomials introduced in \cite{VZ_big}.

\section{Dunkl shift operators on the uniform grid and their polynomial eigensolutions}
\setcounter{equation}{0} In complete analogy with the differential
case \cite{VZ_Bochner}, let us consider the most general linear
operator $L$ of first order with respect to the shift operator
$T^+$  which  also contains the reflection operator $R$: \be L=
F_0(x) + F_1(x) R + G_0(x) T^+ + G_1(x) T^+ R, \lab{gen_L} \ee
where $F_0(x), F_1(x), G_0(x), G_1(x)$ are arbitrary functions.
The shift operator is defined as usual by:
$$
T^+ = \exp(\partial_x)
$$
so that $T^+ f(x) = f(x+1)$ for any function $f(x)$. The operator
$T^+R$ acts on functions $f(x)$ according to $T^+R f(x) =
f(-x-1)$.

We are seeking orthogonal polynomial eigensolutions of the operator $L$, i.e. for every $n$ we assume that there exists
 a monic polynomial $P_n(x)=x^n + O(x^{n-1})$ which is an eigenfunction of the operator $L$ with eigenvalue $\lambda_n$ \re{LPP0}.
 In what follows we will suppose that \be \lambda_n \ne 0  \quad \mbox{for} \;  n=1,2,\dots,
 \quad \lambda_n \ne \lambda_m \quad \mbox{for} \; n \ne m. \lab{ndeg_lambda} \ee
We first establish the necessary conditions for the existence of such eigensolutions.

We will assume the operator $L$  to be potentially self-adjoint.
This leads to the condition \be G_0(x)=0 \lab{G00} \ee (see also
Sect.5 for further details). On the one hand, the formal adjoint
of the operator $T^+$ is the backward-shift operator
$(T^+)^*=T^-$, where
$$
T^- f(x) = f(x-1).
$$
On the other hand, for the adjoint of the operator $T^+ R$ we have
$$
(T^+ R)^* =R^* T^- = R T^- = T^+R .
$$
Thus, the operator $T^+R$ is (formally) self-adjoint, while the
operator $T^+$ is not. Hence the symmetric (or potentially
self-adjoint) operators of type \re{gen_L} cannot contain a term
with the operator $T^+$ only.

Moreover, without loss of generality we can assume that
$F_0+F_1+G_1=0$. Indeed, $L\{1\}=constant$. We can assume that
$L\{1\}=0$ (otherwise we can add a constant to $F_0(x)$ to ensure
this). This leads to the desired condition.

We thus can present the operator $L$ in the form \be L = F(x)
(I-R) + G(x)(T^+ R-I) \lab{L_FG} \ee with only two unknown
functions $F(x), G(x)$. Here $I$ is identity operator.

Equivalently, we can write the action of the operator $L$ on
functions $f(x)$ as: \be L f(x) = F(x)(f(x) - f(-x)) +
G(x)(f(-x-1) -f(x)) . \lab{Lf_def} \ee Considering the action of
the operator $L$ on $x$ and $x^2$ we arrive at the necessary
conditions \be F(x) = \frac{q_1(x) +q_2(x)}{2x}, \quad G(x) =
\frac{q_2(x)}{2x+1}, \lab{FG_expr} \ee where by $q_i(x), \: i=1,2$
we mean arbitrary polynomials of degree $i$, i.e. \be q_1(x) =
\xi_1 x + \xi_0, \quad q_2(x) = \eta_2 x^2 + \eta_1 x + \eta_0
\lab{q1q2} \ee with arbitrary coefficients $\xi_0, \dots, \eta_2$.

It is easily seen that these conditions are also sufficient;
namely, that the operator $L$ defined by \re{L_FG} with functions
$F(x), G(x)$ given by \re{FG_expr}, \re{q1q2}, preserves the
linear space of polynomials. More exactly, for any polynomial
$Q(x)$ of degree $n$ we have $L Q(x) = \t Q(x)$, where $\t Q(x)$
is another polynomial of degree $n$. For monomials $x^n$ we find
\be L x^n = \lambda_n x^n + O(x^{n-1}), \lab{L_x_n} \ee where \be
\lambda_n = \left\{ \frac{\eta_2 n}{2} \quad \mbox{if $n$}  \;
\mbox{is even} \atop      \xi_1 - \frac{\eta_2(n-1)}{2} \quad
\mbox{if $n$}  \; \mbox{is odd} .  \right . \lab{lambda_n_eo} \ee
From this expression it is clear that necessarily \be \eta_2 \ne
0, \quad \xi_1 \ne \eta_2 N, \; N=0,1,2,\dots
\lab{ndeg_par} \ee (otherwise the operator $L$ becomes degenerate:
it does not preserve the degree of the polynomial).

Assuming that conditions \re{ndeg_par} are valid, we can construct
the eigensolutions \be L P_n(x) =\lambda_n P_n(x), \lab{LPP} \ee
for every $n =0,1,2,\dots$, where $P_n(x) =x^n +O(x^{n-1}) $ is a
monic polynomial of degree $n$.

Since $\eta_2 \ne 0$, we can always assume that $\eta_2=1$. We can
thus present the functions $F(x)$ and $G(x)$ as \be G(x) =
\frac{(x-r_1+1/2)(x-r_2+1/2)}{2x+1}, \quad F(x) =
\frac{(x-\rho_1)(x-\rho_2)}{2x} \lab{FG_simpl} \ee with 4
arbitrary real parameters $r_1, r_2,\rho_1, \rho_2$. We will use
this parametrization in what follows.

The eigenvalue $\lambda_n$ has the expression \be \lambda_n =
\left\{ \frac{n}{2} \quad \mbox{if $n$}  \; \mbox{is even} \atop
r_1+r_2-\rho_1-\rho_2 - \frac{(n+1)}{2} \quad \mbox{if $n$}  \;
\mbox{is odd} .  \right . \lab{lambda_fix} \ee From the properties
of the operator $L$ described above, it follows that the monic
polynomials $P_n(x)$ exist for all $n=0,1,2,\dots$ and are the
unique polynomial solutions of equation \re{LPP} for the given
eigenvalues $\lambda_n$.

The polynomials $P_n(x)=P_n(x;r_1,r_2, \rho_1,\rho_2)$ depend on 4
parameters $r_1,r_2, \rho_1,\rho_2$. There is however an obvious
invariance under the action of the Klein group $Z_2 \times Z_2$.
Indeed, the permutations $r_1 \leftrightarrow r_2$ and $\rho_1
\leftrightarrow \rho_2$ leave the operator $L$ unchanged and hence
the polynomials $P_n(x;r_1,r_2, \rho_1,\rho_2)$ are invariant as
well under these permutations. (With the identity this $Z_2 \times
Z_2$ group contains 4 elements.)

\section{A 2-diagonal basis and an explicit expression in terms of hypergeometric functions}
\setcounter{equation}{0} Choose the following polynomial basis.
For even degree, take \be \phi_{2n}(x)= (x-r_1+1/2)_n(-x-r_1+1/2)_n
\lab{phi_ev} \ee and for odd degree, take \be \phi_{2n+1}(x) =
(x-r_1+1/2)_{n+1}(-x-r_1+1/2)_n =
(x-r_1+1/2)(x-r_1+3/2)_{n}(-x-r_1+1/2)_n. \lab{phi_odd} \ee Here
$(x)_n=x(x+1)(x+2) \dots (x+n-1)$ stands as usual for the shifted
factorial (Pochhammer symbol). Obviously, $\phi_n(x)$ is a
polynomial of degree $n$: \be \phi_n(x) =(-1)^{n(n-1)/2} x^n +
O(x^{n-1}). \lab{psi_leading} \ee
 It is directly verified that the operator $L$ is two-diagonal in this basis
\be L \phi_n(x) = \lambda_n \phi_n(x) + \nu_n \phi_{n-1}(x),
\lab{L_2_diag} \ee where \be \nu_n = \left\{ \frac{n}{2}(r_1+r_2-
\frac{n}{2}) \quad \mbox{if $n$}  \; \mbox{is even} \atop
(\rho_1-r_1+n/2)(\rho_2-r_1+n/2)  \quad \mbox{if $n$}  \; \mbox{is
odd}   \right . . \lab{nu_n} \ee Using this striking observation,
we can explicitly construct  the polynomial eigensolutions of the
operator $L$ in the same manner as in \cite{VZ_big}.

Let us expand the polynomials $P_n(x)$ over the basis $\phi_n(x)$:
$$
P_n(x) = \sum_{s=0}^n A_{ns} \phi_s(x).
$$
From the eigenvalue problem \re{LPP}, we have the recurrence
relation for the coefficients $A_{ns}$: \be A_{n,s+1} =
\frac{A_{ns}(\lambda_n-\lambda_s)}{\nu_{s+1}}. \lab{rec_A_s} \ee
Hence the coefficients $A_{ns}$ can easily be found in terms of
$A_{n0}$: \be A_{ns}= A_{n0} \:
\frac{(\lambda_n-\lambda_0)(\lambda_n-\lambda_1) \dots
(\lambda_n-\lambda_{s-1})}{\nu_1 \nu_2 \dots \nu_s} \lab{A_ns_A_0}
\ee or in terms of the coefficient $A_{nn}$: \be A_{ns}= A_{nn} \:
\frac{\nu_n \nu_{n-1} \dots
\nu_{s+1}}{(\lambda_n-\lambda_{n-1})(\lambda_n - \lambda_{n-2})
\dots (\lambda_n - \lambda_s) }. \lab{A_ns_A_n} \ee We thus have
the following explicit formula for the polynomials $P_n(x)$ \be
P_n(x) = A_{n0} \: \sum_{s=0}^n
\frac{(\lambda_n-\lambda_0)(\lambda_n-\lambda_1) \dots
(\lambda_n-\lambda_{s-1})}{\nu_1 \nu_2 \dots \nu_s} \phi_s(x).
\lab{P_n_hyp_phi} \ee Expression \re{P_n_hyp_phi} resembles Gauss'
hypergeometric function and can be considered as a nontrivial
generalization of it.

The expansion coefficients $A_{ns}$ have different expressions
depending on the parity of the numbers $n$ and $s$.

When $n$  is even $n=2,4,6,\dots$ we have \be A_{n,2s}/A_{n0} =
\frac{(-n/2)_s (n/2+\rho_1+\rho_2-r_1-r_2)_s}{s!
(1-r_2-r_1)_s(1/2+\rho_1-r_1)_s(1/2+\rho_2-r_1)_s}, \quad s=0,1,
\dots n/2 \lab{A_ee} \ee and \be A_{n,2s+1}/A_{n0} = \xi_n \;
\frac{(1-n/2)_s   (n/2+\rho_1+\rho_2-r_1-r_2)_s}{s!
(1-r_2-r_1)_s(3/2+\rho_1-r_1)_s(3/2+\rho_2-r_1)_s}, \quad s=0,1,
\dots n/2-1, \lab{A_eo} \ee where
$$
\xi_n=\frac{n}{2(1/2+\rho_1-r_1)(1/2+\rho_2-r_1)}.
$$
When $n$ is odd $n=1,3,5,\dots$  then \be A_{n,2s}/A_{n0} =
\frac{((1-n)/2)_s ((n+1)/2+\rho_1+\rho_2-r_1-r_2)_s}{s!
(1-r_2-r_1)_s(1/2+\rho_1-r_1)_s(1/2+\rho_2-r_1)_s}, \quad s=0,1,
\dots (n-1)/2 \lab{A_oe} \ee and \be A_{n,2s+1}/A_{n0} = -\eta_n
\; \frac{((1-n)/2)_s   ((n+3)/2+\rho_1+\rho_2-r_1-r_2)_s}{s!
(1-r_2-r_1)_s(3/2+\rho_1-r_1)_s(3/2+\rho_2-r_1)_s}, \quad s=0,1,
\dots (n-1)/2, \lab{A_oo} \ee where
$$
\eta_n = \frac{(n+1)/2
+\rho_1+\rho_2-r_1-r_2}{(1/2+\rho_1-r_1)(1/2+\rho_2-r_1)}.
$$
It is now convenient to separate the even and odd terms in the
expression for the polynomials: \be P_n(x) = \sum_{s=0} A_{n,2s}
\phi_{2s}(x) + \sum_{s=0} A_{n,2s+1} \phi_{2s+1}(x). \lab{P_evod}
\ee Taking into account the explicit formulas \re{A_ee}-\re{A_oo}
for the expansion coefficients, we finally obtain the following
expressions in terms of hypergeometric functions:

(i) if $n$ is even \ba &&\frac{P_n(x)}{A_{n0}}= {_4}F_3 \left(
{-n/2, n/2+1+\rho_1+\rho_2-r_1-r_2, x-r_1+1/2,-x+r_1+1/2 \atop
1-r_1-r_2, 1/2+\rho_1-r_1, 1/2+\rho_2-r_1 } ; 1 \right) +
\nonumber \\  \lab{hyp_P_ev}
\\&&\xi_n (x-r_1+1/2) \; {_4}F_3 \left(  {1-n/2,
n/2+1+\rho_1+\rho_2-r_1-r_2, x-r_1+3/2,-x+r_1+1/2 \atop 1-r_1-r_2,
3/2+\rho_1-r_1, 3/2+\rho_2-r_1 } ; 1 \right); \nonumber \ea

(ii) if $n$ is odd \ba &&\frac{P_n(x)}{A_{n0}}= {_4}F_3 \left(
{-(n-1)/2, (n+1)/2+\rho_1+\rho_2-r_1-r_2, x-r_1+1/2,-x+r_1+1/2
\atop  1-r_1-r_2, 1/2+\rho_1-r_1, 1/2+\rho_2-r_1 } ; 1 \right) -
\nonumber \\  \lab{hyp_P_od}
\\&&\eta_n\: (x-r_1+1/2) \; {_4}F_3 \left(  {-(n-1)/2,
(n+3)/2+\rho_1+\rho_2-r_1-r_2, x-r_1+3/2,-x+r_1+1/2 \atop
1-r_1-r_2, 3/2+\rho_1-r_1, 3/2+\rho_2-r_1 } ; 1 \right) .
\nonumber \ea

Note that we have terminating hypergeometric functions (i.e. only
a finite number of terms in the hypergeometric series appear).
Hence, as is easily seen, the expressions in the rhs of
\re{hyp_P_ev} and \re{hyp_P_od} are polynomials of degree $n$ in
the argument $x$. The coefficient $A_{n0}$ (ensuring that the
polynomials $P_n(x)$ are monic) can be found from the coefficient
in front of the leading term $x^n$.

{\it Remark .} At first sight one might be tempted to identify the
hypergeometric functions ${_4}F_3(1)$ with the Wilson (Racah)
polynomials  which are known to be expressible in terms of
truncated ${_4}F_3(1)$ functions \cite{KLS}. This is not so,
however, because the balance is not right.

Recall that the generalized hypergeometric function
$$
{_p}F_q \left(  { a_1, a_2, \dots a_p \atop b_1, b_2, \dots b_q};1 \right)
$$
is said to be $k$-balanced if $a_1+a_2 + \dots + a_p = k + b_1 +
b_2 + \dots + b_q$.  The Racah polynomials are expressed in terms
of -1 balanced hypergeometric functions ${_4}F_3(1)$ \cite{KLS}.
In \re{hyp_P_ev}, \re{hyp_P_od}, as is easily seen, all
hypergeometric functions ${_4}F_3(1)$  are instead zero-balanced.
In what follows we will give another explicit formula for the
polynomials $P_n(x)$ which will relate them to "true" Wilson
polynomials.

Remark that the basis $\phi_n(x)$ is important in Terwilliger's
approach to Leonard duality \cite{Ter}, \cite{Ter2}. Similar bases
exist for -1 Jacobi polynomials as well \cite{VZ_big}.

Equivalently, the basis $\phi_n(x)$ can be presented as $$
\phi_n(x) = (-1)^{n(n-1)/2} \prod_{i=1}^n (x -\alpha_i), $$ where
$$
\alpha_n=(-1)^n \left(r_1 + \frac{n}{2}-\frac{1}{4}
\right)-\frac{1}{4}.
$$
This means that the $\phi_n(x)$ form a special case of the
Newtonian interpolation polynomial basis  $\Omega_n(x) =
\prod_{i=1}^n (x -\alpha_i)$, where $\alpha_i$ are fixed
interpolation nodes \cite{Hild}. In the theory of orthogonal
polynomials such bases were first proposed by Geronimus
\cite{Ger2}. All known explicit examples of orthogonal polynomials
(including the Askey-Wilson polynomials at the top level) have
simple expansion coefficients with respect to a Newtonian basis
$\Omega_n(x)$ with appropriately chosen nodes $\alpha_i$
\cite{Ter}, \cite{Ter2}. We thus see that the BI polynomials
satisfy this property as well.

\section{Bannai-Ito algebra, structure relations and recurrence coefficients for the BI polynomials}
\setcounter{equation}{0} In this section we derive the 3-term
recurrence relations as well as the structure relations for the
BI-polynomials. The main tool will be the Bannai-Ito algebra
(BI-algebra for brevity) which is a special case of the AW(3)
algebra introduced in \cite{Zhe}. For the Askey-Wilson polynomials
the structure relations could be derived in a similar manner using
representations of the AW(3) algebra \cite{K_Str},
\cite{Ismail_Str}. For a study of the structure relations of the
-1 Jacobi polynomials see \cite{TVZ_J}.

Consider the operators \be X= 2L + \kappa, \quad Y=2x+1/2,
\lab{def_XY} \ee where by $x$ we mean the operator multiplication
by $x$ and where
$$
\kappa=\rho_1+\rho_2-r_1-r_2+1/2.
$$
We define also the third operator \be Z= \frac{r_1r_2}{x+1/2} +
\frac{\rho_1 \rho_2}{x} - \frac{(x-\rho_1)(x-\rho_2)}{x} R -
\frac{(x-r_1+1/2)(x-r_2+1/2)}{x+1/2} T^+R . \lab{Z_def} \ee It is
easily verified that the operator $Z$ transforms any polynomial of
degree $n$ into a polynomial of degree $n+1$ and that \be Z x^n =
2(-1)^{n+1}  x^{n+1} + O(x^n) . \lab{Zx^n} \ee

The operators X, Y, Z are seen to have the following Jordan or
anticommutator products which are taken to be the defining
relations of the BI algebra: \be \{X,Y\}=Z+\omega_3,  \quad
\{Z,Y\}=X+\omega_1, \quad \{X,Z\}=Y+\omega_2, \lab{XYZ_BI} \ee
where
$$
\omega_1=4(\rho_1\rho_2+ r_1 r_2), \quad \omega_2= 2(\rho_1^2
+\rho_2^2  -r_1^2-r_2^2), \quad  \omega_3=4 (\rho_1\rho_2-
r_1r_2).
$$
(Note that the first relation \re{XYZ_BI} can be considered as a definition of the operator $Z$).

In the simplest case $\omega_1 = \omega_2=\omega_3=0$, the BI
algebra was considered in \cite{AK}, \cite{GP}, \cite{OS},  \cite{RW} as an anticommutator analogue
of the rotation algebra $so(3)$. A contraction of \re{XYZ_BI} has
also been used in \cite{TVZ_J} to provide an algebraic account of
structural properties of the big -1 Jacobi polynomials.

The Casimir operator $Q$ commuting with all generators $X,Y,Z$ of
the BI algebra has a simple form: \be Q=X^2+Y^2+Z^2. \lab{Casimir}
\ee In the given realization of the operators $X,Y,Z$ in terms of
difference-reflection operators, this Casimir operator takes the
value \be Q= 2(\rho_1^2 + \rho_2^2 + r_1^2+r_2^2) -1/4 .
\lab{Q_value} \ee

The BI algebra can be exploited to find explicitly the 3-term
recurrence relation of the polynomials $P_n(x)$.

We first introduce  two important operators $J_{+}$ and $J_{-}$ by
the formulas \be J_{+}=(Y+Z)(X-1/2)-\frac{\omega_2+\omega_3}{2}
\lab{J+} \ee and \be
J_{-}=(Y-Z)(X+1/2)+\frac{\omega_2-\omega_3}{2} . \lab{J-} \ee From
the  commutation relations \re{XYZ_BI}, we find that the operators
$J_{\pm}$ satisfy the (anti)commutation relations \be \{X,J_{+}\}
= J_{+}, \quad \{X,J_{-}\} = -J_{-}. \lab{com_JPM} \ee Moreover,
from \re{com_JPM} it is seen that both $J_{+}^2$ and $J_{-}^2$
commute with the operator $X$: \be [X,J_{+}^2]=[X,J_{-}^2]=0.
\lab{com_X_JPM} \ee The operator $J_{+}$ has the following
property: it annihilates any constant $J_{+} \{1\}=0$ and
transforms any monomial $x^n$ of even degree n into a polynomial
of the same degree $n$ (providing that $r_1+r_2 \ne n/2$): \be
J_{+} x^n = n(n-2r_1-2r_2) x^{n} + O(x^{n-1}), \quad n=0,2,4,
\dots, \lab{J+_even} \ee while any monomial $x^n$ of odd degree is
transformed into a polynomial of even degree $n+1$ (providing that
$r_1+r_2-\rho_1-\rho_2 \ne n+1$) \be J_{+} x^n =
4(r_1+r_2-\rho_1-\rho_2-n-1) x^{n+1} + O(x^{n}), \quad
n=1,3,5,\dots .  \lab{J+_odd} \ee

Quite similarly, the operator $J_{-}$ transforms any monomial
$x^n$ of even degree $n$ into a polynomial of odd degree $n+1$:
\be J_{-} x^n = 4(\rho_1+\rho_2-r_1-r_2+n+1) x^{n+1} + O(x^{n}),
\quad n=0,2,4, \dots, \lab{J-_even} \ee while the monomial $x^n$
of odd degree is transformed to a polynomial of the same degree
$n$: \be J_{-} x^n = \left(4(\rho_1\rho_2-r_1r_2) -n^2 +
2n(r_1+r_2) \right) x^{n} + O(x^{n-1}) , \quad n=1,3,5,\dots .
\lab{J-_odd} \ee A careful analysis (where the leading and next to
leading monomials are considered), shows that the operator $J_{+}$
transforms any polynomial of even degree $n$ into a polynomial of
degree $n$ or less, while it maps any polynomial of odd degree
into a polynomial of exact degree $n+1$; similarly, it is seen
that the operator $J_{-}$ transforms any polynomial of even degree
$n$ into a polynomial of exact degree $n+1$, while it maps any
polynomial of odd degree into a polynomial of degree $n$ or less.

Now, let $\psi_n(x)$ be any eigenfunction of the operator $X$:
$$
X \psi_n(x) = \mu_n \psi_n(x) ,
$$
where \be \mu_n= 2\lambda_n +\kappa =(-1)^n (n+\kappa). \lab{mu}
\ee  It then follows from relations \re{com_JPM} that the function
$\t \psi_n(x)=J_{+}\psi_n$ will again be an eigenfunction of the
operator $X$ corresponding to the eigenvalue $\t \mu_n = 1-\mu_n$.
It is seen that $\t \mu_n = \mu_{n-1}$ if $n$ is even and that $\t
\mu_n = \mu_{n+1}$ if $n$ is odd.

A similar property is found for the operator $J_{-}$. In this case
the function $J_{-} \psi_n$ will be an eigenfunction of the
operator $X$ with eigenvalue $-1-\mu_n$ which is $\mu_{n+1}$ for
even $n$ and $\mu_{n-1}$ for odd $n$.

We have already established that the operators $J_{\pm}$ transform
polynomials into polynomials. Hence from the above observations,
\be J_{+}P_n(x) = \left\{ {\alpha_n^{(0)} P_{n-1}(x), \quad
\mbox{if} \quad n \quad \mbox{even},  \atop \alpha_n^{(1)}
P_{n+1}(x), \quad \mbox{if} \quad n \quad \mbox{odd}} \right .
\lab{J+_P} \ee and similarly, \be J_{-}P_n(x) = \left\{
{\beta_n^{(0)} P_{n+1}(x), \quad \mbox{if} \quad n \quad
\mbox{even}  \atop \beta_n^{(1)} P_{n-1}(x), \quad \mbox{if} \quad
n \quad \mbox{odd}} \right . . \lab{J-_P} \ee The coefficients in
the rhs of these formulas can be found by comparison of
highest-order terms:
$$
\alpha_n^{(0)} = \frac{2n(\rho_1+\rho_2
+n/2)(r_1+r_2-n/2)(r_1+r_2-\rho_1-\rho_2-n/2)}{r_1+r_2-\rho_1-\rho_2-n}
\quad, \alpha_n^{(1)}=4(r_1+r_2-\rho_1-\rho_2-n-1)
$$
and
$$
\beta_n^{(0)}=-4(r_1+r_2-\rho_1-\rho_2-n-1)=-\alpha_n^{(1)}, \quad
\beta_n^{(1)}=
\frac{4(\rho_1-r_1+n/2)(\rho_1-r_2+n/2)(\rho_2-r_1+n/2)(\rho_2-r_2+n/2)}{\rho_1+\rho_2-r_1-r_2+n}
.
$$

Note that formulas \re{J+_P} and \re{J-_P} show that the operators
$J_{\pm}$ are block-diagonal in the basis of the  polynomials
$P_n(x)$ with each block a $2\times 2$ matrix.

From these formulas it follows that the operators $J_{\pm}^2$ have
the polynomials $P_n(x)$ as eigenfunctions \be J_{+}^2 P_n(x)
=\left\{ {\alpha_n^{(0)} \alpha_{n-1}^{(1)} P_{n}(x), \quad
\mbox{if} \quad n \quad \mbox{even}  \atop \alpha_{n+1}^{(0)}
\alpha_n^{(1)} P_{n}(x), \quad \mbox{if} \quad n \quad \mbox{odd}}
\right . \lab{J+2_P} \ee and similarly \be J_{-}^2 P_n(x) =\left\{
{\beta_n^{(0)} \beta_{n+1}^{(1)} P_{n}(x), \quad \mbox{if} \quad n
\quad \mbox{even}  \atop \beta_{n-1}^{(0)} \beta_n^{(1)} P_{n}(x),
\quad \mbox{if} \quad n \quad \mbox{odd}} \right . . \lab{J-2_P}
\ee Now, from the commutation relations \re{XYZ_BI} we find \be
J_{+}^2=(Y+Z)^2(-X^2+X-1/4) +(\omega_2+\omega_3)^2/4 .
\lab{J+2_rel1} \ee Using formulas  \re{Casimir} and \re{XYZ_BI} we
also see that
$$
(Y+Z)^2 = Y^2+Z^2 +\{Y,Z\} = Q-X^2+X+\omega_1 .
$$
We thus have \be J_{+}^2=(Q-X^2+X+\omega_1)(-X^2+X-1/4)
+(\omega_2+\omega_3)^2/4, \lab{J+2_rel2} \ee where $Q$ is the
Casimir operator \re{Casimir}. Taking into account that in our
realization of the BI algebra the Casimir operator takes the value
\re{Q_value}, we arrive at the expression \be
J_{+}^2=\left((X+\rho_1+\rho_2-1/2)^2 -(r_1+r_2)^2  \right)
\left((X-\rho_1-\rho_2-1/2)^2 -(r_1+r_2)^2  \right) . \lab{J+2_4p}
\ee Quite similarly we obtain \be
J_{-}^2=\left((X+\rho_2-\rho_1+1/2)^2 -(r_2-r_1)^2  \right)
\left((X+\rho_1-\rho_2+1/2)^2 -(r_2-r_1)^2  \right) . \lab{J-2_4p}
\ee We see that both operators $J_{+}^2$ and $J_{-}^2$ (which
commute with the operator $X$) are  expressible as 4-th degree
polynomials in $X$.

Starting from relations \re{J+_P} and \re{J-_P} one can now derive
many useful relations for the polynomials $P_n(x)$.

As a first example we present the structure relations of the
polynomials $P_n(x)$.

Define the operator $U_n^{(1)}, \: n=0,1,2,\dots$ as \be U_n^{(1)}
= \left\{ {J_{+}, \quad \mbox{if} \quad n \quad \mbox{even}  \atop
J_{-}, \quad \mbox{if} \quad n \quad \mbox{odd}} \right . .
\lab{U1} \ee
 We have in this case
 \be
 U_n^{(1)} P_n(x) =\epsilon^{(1)}_n P_{n-1}(x), \quad \epsilon^{(1)}_n =\left\{{\alpha_n^{(0)} , \quad \mbox{if} \quad n \quad
 \mbox{even} \atop  \beta_n^{(1)}, \quad \mbox{if} \quad n \quad \mbox{odd}} \right . . \lab{U1P} \ee
Similarly, define the operator $U_n^{(2)}, \: n=0,1,2,\dots$ as
\be U_n^{(2)} = \left\{ {J_{-}, \quad \mbox{if} \quad n \quad
\mbox{even}  \atop J_{+}, \quad \mbox{if} \quad n \quad
\mbox{odd}} \right . . \lab{U2} \ee
 Then we have
 \be
 U_n^{(2)} P_n(x) =\epsilon^{(2)}_n P_{n+1}(x), \quad \epsilon^{(1)}_n =\left\{{\beta_n^{(0)} , \quad \mbox{if} \quad n \quad
 \mbox{even} \atop  \alpha_n^{(1)}, \quad \mbox{if} \quad n \quad \mbox{odd}} \right . . \lab{U2P} \ee
It is seen that the operator $U_n^{(1)}$ plays the role of the
lowering operator, while the operator  $U_n^{(2)}$ serves as the
raising operator for the polynomials $P_n(x)$.

We are now ready to derive, as a second example, the 3-term
recurrence relations of the polynomials $P_n(x)$.

To do this we consider the operator \be
V=J_{+}(X+1/2)+J_{-}(X-1/2) =2Y(X^2-1/4) -\omega_3 X-\omega_2/2 .
\lab{V_def} \ee By construction, this operator is 2-diagonal in
the basis of the polynomials $P_n(x)$: \be V P_n(x) =\left\{{
(\mu_n+1/2)\alpha_n^{(0)}P_{n-1}(x)
+(\mu_n-1/2)\beta_n^{(0)}P_{n+1}(x), \quad \mbox{if} \quad n \quad
\mbox{even} \atop   (\mu_n-1/2)\beta_n^{(1)}P_{n-1}(x)
+(\mu_n+1/2)\alpha_n^{(1)}P_{n+1}(x) \quad \mbox{if} \quad n \quad
\mbox{odd}} \right . , \lab{V_2term} \ee where $\mu_n$ is the
eigenvalue of the operator $X$ given in \re{mu}.

Moreover, from \re{def_XY}, the operator $Y$ coincides (up to an
affine transformation) with the multiplication by $x$, that is, in
the polynomial basis $P_n(x)$, we have $Y P_n(x) = (2x
+1/2)P_n(x)$. Hence, by the second equality in \re{V_def}, \be V
P_n(x) = ((\mu_n^2-1/4)(4x+1) -\omega_3 \mu_n -\omega_2/2)P_n(x).
\lab{VP_x} \ee Comparing \re{V_2term} and \re{VP_x}, we arrive at
the 3-term recurrence relation for the polynomials $P_n(x)$: \be x
P_n(x) = P_{n+1}(x) + b_n P_n(x) + u_n P_{n-1}(x), \lab{3-term_P}
\ee where the coefficients are \be b_n = -1/4 + \frac{\omega_3
\mu_n + \omega_2/2}{4\mu_n^2-1} \lab{bBI} \ee and \be
u_n=\left\{{\frac{\alpha_n^{(0)}}{4 \mu_n -2} , \quad \mbox{if}
\quad n \quad         \mbox{even} \atop \frac{\beta_n^{(1)}}{4
\mu_n +2}  , \quad \mbox{if} \quad n \quad \mbox{odd}} \right . .
\lab{uBI} \ee We can also present the expressions for the
coefficients $u_n$ in a more detailed form. For even
$n=0,2,4,\dots$ we have \be u_n=-{\frac {n \left(
n+2\,\rho_{{1}}+2\,\rho_{{2}} \right) \left(
n-2\,r_{{1}}-2\,r_{{2}} \right)  \left(
n+2\,\rho_{{1}}+2\,\rho_{{2}}- 2\,r_{{1}}-2\,r_{{2}} \right) }{
16\left( n+\rho_{{1}}+\rho_{{2}}-r_{{1} }-r_{{2}} \right) ^{2}}}.
\lab{even_u} \ee
 For odd $n=1,3,5,\dots$ we have
\be
u_n = -{\frac { \left( n+2\,\rho_{{1}}-2\,r_{{1}} \right)  \left( n+2
\,\rho_{{1}}-2\,r_{{2}} \right)  \left( n+2\,\rho_{{2}}-2\,r_{{1}}
 \right)  \left( n+2\,\rho_{{2}}-2\,r_{{2}} \right) }{ 16\left( n+\rho_{
{1}}+\rho_{{2}}-r_{{1}}-r_{{2}} \right) ^{2}}}. \lab{odd_u} \ee We
can now identify these recurrence coefficients with those of the
Bannai-Ito polynomials.

Recall that the monic Bannai-Ito polynomials can be defined
through the 3-term recurrence relation \re{3-term_P}, with the
recurrence coefficients \cite{BI}, \cite{VZ_big} \be u_n = A_{n-1}
C_n, \quad b_n = \theta_0 -A_n-C_n, \lab{ub_AC} \ee where \be
A_n=\left\{ { \frac{2h(n+1+\tau_1)(n+1+\tau_2)}{2n+2-s^*}, \quad
\mbox{} \quad n \quad \mbox{even} \atop
\frac{2h(n+1-s^*)(n+1-\tau_3)}{2n+2-s^*}  \quad \mbox{} \quad n
\quad \mbox{odd}} \right . \lab{A_BI} \ee and \be C_n=\left\{ {
-\frac{2h n(n-s^*+\tau_3)}{2n-s^*}, \quad \mbox{} \quad n \quad
\mbox{even} \atop -\frac{2h(n-\tau_1-s^*)(n-\tau_2-s^*)}{2n-s^*}
\quad \mbox{} \quad n \quad \mbox{odd}} \right . .\lab{C_BI} \ee
In these formulas the parameters $\tau_1, \tau_2, \tau_3, s^*$ are
the 4 essential parameters of the Bannai-Ito polynomials while the
parameters $h$ and $\theta_0$ define an affine transformation $x
\to (x + \theta_0)h^{-1}$ of the argument of the polynomials
$P_n(x)$.

An easy analysis shows that the coefficients \re{bBI} and \re{uBI}
coincide with the recurrence coefficients \re{ub_AC} under the
identifications \be \theta_0 =\rho_1, \: h=1/4, \:
\tau_1=2(\rho_1-r_1), \: \tau_2=2(\rho_1-r_2), \:
\tau_3=-2(\rho_1+\rho_2), \: s^*=2(r_1+r_2-\rho_1-\rho_2).
\lab{param_BI} \ee Note that apart from the solution
\re{param_BI}, there exist 3 other solutions which can be obtained
from \re{param_BI} by the permutations $r_1 \leftrightarrow r_2,
\; \rho_1 \leftrightarrow \rho_2$  . This corresponds to the
previously mentioned symmetry of the BI polynomials with respect
to the $Z_2 \times Z_2$ group.

We thus see that the orthogonal polynomials $P_n(x)$ coincide with
the general Bannai-Ito polynomials.

For the reader's convenience we record the explicit formulas for
the recurrence coefficients: \be P_{n+1}(x) + (\rho_1 - A_n -
C_n)P_n(x) + A_{n-1}C_n P_{n-1}(x) =x P_n(x), \lab{rec_AC_BIP} \ee
where \be A_n=\left\{ { {\frac { \left( n+1+2\,{\rho_1}-2\,{r_1}
\right) \left( n+1 +2\,{\rho_1}-2\,{ r_2} \right) }{4(n+1-{ r_1}-{
r_2}+{\rho_1}+ {\rho_2})}} , \quad \mbox{} \quad n \quad
\mbox{even} \atop {\frac { \left( n+1-2\,{r_1}-2\,{
r_2}+2\,{\rho_1}+2\,{ \rho_2} \right) \left(
n+1+2\,{\rho_1}+2\,{\rho_2} \right) }{4(n+
1-{r_1}-{r_2}+{\rho_1}+{\rho_2})}}
  ,\quad\mbox{} \quad n \quad
\mbox{odd}} \right . \lab{A_BIP} \ee and \be C_n=\left\{ {
-{\frac {n \left( n-2\,{r_1}-2\,{r_2} \right) }{4(n-{r_1}-
{r_2}+{\rho_1}+{\rho_2})}}
, \quad \mbox{} \quad n \quad
\mbox{even} \atop -{\frac { \left( n-2\,{r_2}+2\,{\rho_2} \right)  \left( n-2
\,{r_1}+2\,{\rho_2} \right) }{4(n-{r_1}-{r_2}+{\rho_1}+{
\rho_2})}},
 \quad
\mbox{} \quad n \quad \mbox{odd}.} \right . . \lab{C_BIP} \ee

Consider now the question of the positive definiteness of the
polynomials $P_n(x)$ obeying a recurrence relation of the form
\re{3-term_P}. Recall that polynomials $P_n(x)$ are called
positive definite if $u_n>0$ for all $n=1,2,\dots$. This property
is equivalent to the existence of a positive measure $d \mu(x)$ on
the real axis such that the orthogonality property reads \be
\int_a^b P_n(x) P_m(x) d \mu(x) = h_n, \delta_{nm} \lab{ort_P} \ee
where $h_n =u_1 u_2 \dots u_n$, and  where the integration limits
$a,b$ can be either finite or infinite.

In our case it is seen that $u_n<0$ for sufficiently large $n$. It
is hence impossible to ensure positivity for all $n=1,2,\dots$.
Nevertheless it is possible to obtain a finite set of positive
definite orthogonal polynomials.

If $u_i>0$ for $i=1,2,\dots, N$ and  $u_{N+1}=0$, it is well known
that we have a finite system of orthogonal polynomials $P_0(x),
P_1(x), \dots, P_N(x)$ satisfying the discrete orthogonality
relation \be \sum_{s=0}^N w_s P_n(x_s) P_m(x_s) = h_n \:
\delta_{nm}, \quad h_n = u_1 u_2 \dots u_n, \lab{N-ort} \ee where
$x_s, \: s=0,1,\dots,N$ are the simple roots of the polynomial
$P_{N+1}(x)$. (The fact that the roots $x_s$ are simple is
guaranteed by the condition $u_n>0,\: n=1,2,\dots, N$ \cite{Chi}.)

The discrete weights $w_s$ are then given by the formula
\cite{Chi} \be w_s = \frac{P_{N}^{(1)}(x_s)}{P_{N+1}'(x_s)} =
\frac{h_N}{P_{N}(x_s) P_{N+1}'(x_s)}>0, \quad s=0,1,\dots,N,
\lab{w_s_def} \ee where $P_n^{(1)}(x)$ are associative polynomials
satisfying the recurrence relation
$$
P_{n+1}^{(1)}(x) + b_{n+1} P_{n}^{(1)}(x)+ u_{n+1}
P_{n-1}^{(1)}(x)=x P_{n}^{(1)}(x)
$$
with initial conditions $P_{0}^{(1)}=1, \; P_{1}^{(1)}(x)=x-b_1$.

We thus have positive definite polynomials $P_n(x)$ which are
orthogonal on the finite set of distinct points $x_0, x_1, \dots
x_N$.

Let us consider when such a situation takes place for our
polynomials.

First, assume that $N$ is even.  We have then from \re{odd_u} that
the condition $u_{N+1}=0$ is equivalent to one of 4 possible
conditions
$$
2(r_i-\rho_k)=N+1, \quad i,k=1,2.
$$
We restrict ourselves with the condition \be 2(r_2-\rho_2)=N+1 .
\lab{ev_N_rr} \ee Then it is sufficient to take \be
r_2=r_{{1}}+e+\frac{N}{2}, \; \rho_1=r_{{1}}+e+d+\frac{N-1}{2}, \;
\rho_2=r_{{1}}-1/2+e, \lab{cond_rr_ev} \ee where $r_1, e,d,$ are
arbitrary positive parameters. Assuming \re{cond_rr_ev}, we can
present the recurrence coefficients $u_n$ in the following form.
For even $n=2,4,\dots,N$
$$
u_n= {\frac {n \left( n+4\,r_{{1}}-2+4\,e+2\,d+N \right)
 \left( -n+4\,r_{{1}}+2\,e+N \right)  \left( n-2+2\,e+2\,d
 \right) }{16 \left( n-1+e+d \right) ^{2}}},
$$
for odd $n=1,3,\dots, N-1$
$$
u_n={\frac { \left( n-1+2\,e+2\,d+N \right)  \left( n-1+2\,
d \right)  \left( n-1+2\,e \right)  \left( -n+1+N \right) }{16
 \left( n-1+e+d \right) ^{2}}}.
$$
From these expressions it is seen indeed that the positivity
condition $u_n >0$ is satisfied for $n=1,2,\dots, N$.

Consider now the case when $N>1$ is odd. Then from \re{even_u}, we
see that the condition $u_{N+1}=0$ is equivalent to one of three
conditions:

(i) $\rho_1 + \rho_2 =-(N+1)/2,$

(ii) $r_1 + r_2 =(N+1)/2,$

(iii) $\rho_1 + \rho_2 -r_1-r_2 =-(N+1)/2.$

Condition (iii) leads to singular $u_n$ for $n=(N+1)/2$. Hence
only conditions (i) and (ii) are admissible.

Consider condition (ii) for definiteness.

It is convenient to introduce the parametrization
$$
r_1=(1-\zeta-\eta)/2, \; r_2 = (\eta+\zeta+N)/2, \; \rho_1=(\eta-\zeta)/2, \; \rho_2=(\zeta-\eta-2\xi -N+1)/2,
$$
where $\xi, \eta, \zeta$ are positive parameters with the restriction $\xi > \eta$.
Condition (ii) holds automatically and the recurrence coefficients become
$$
u_n={\frac {n \left( -n+N-1+2\,\xi \right)  \left( -n+N+1 \right)
 \left( -n+2\,N+2\,\xi \right) }{ 16 \: \left( -n+N+\xi \right) ^{2}}}
$$
for even $n$ and
$$
u_n={\frac { \left( n+2\,\eta-1 \right)  \left( -n+2\,\zeta+N
 \right)  \left( -n-2\,\zeta+N+2\,\xi \right)  \left( -n+2\,\eta+2\,N-
1+2\,\xi \right) }{16\: \left( -n+N+\xi \right) ^{2}}}
$$
for odd $n$. Again it is visible that $u_n>0$ when $n=1,2,\dots,N$
and $u_{N+1}=0$.

\section{Complementary Bannai-Ito polynomials and an alternative expression in terms of Racah polynomials}
\setcounter{equation}{0} Let us start with the following simple
lemma
\begin{lem}
Let $P_n(x)$ be monic orthogonal polynomials satisfying the
recurrence relation \be P_{n+1}(x) +(\theta - A_n - C_n)P_n(x) +
C_n A_{n-1} P_{n-1}(x) = xP_n(x), \lab{3-term_AC} \ee with the
standard initial conditions
$$
P_0=1, \; P_1(x) = x-\theta+A_0 .
$$
Assume that the real coefficients $A_n,C_n$ are such that
$A_{n-1}>0,\; C_n
>0, \; n=1,2\dots$ and that $C_0=0$. Take $\theta$ to be an arbitrary real
parameter.

Define the new polynomials \be \t P_n(x) = \frac{P_{n+1}(x) - A_n
P_n(x)}{x-\theta}. \lab{tP_CT} \ee Then the monic polynomials $\t
P_n(x)$ are orthogonal and satisfy the recurrence relation \be \t
P_{n+1}(x) +(\theta - A_n - C_{n+1}) \t P_n(x) + C_n A_{n} \t
P_{n-1}(x) = x \t P_n(x). \lab{3-term_tP} \ee The inverse
transformation from the polynomials $\t P_n(x)$ to the polynomials
$P_n(x)$ is given by the formula \be P_n(x) = \t P_n(x) - C_n \t
P_{n-1}(x) . \lab{GT} \ee
\end{lem}
To the best of our knowledge this Lemma is due to Karlin and
McGregor \cite{KMG}.

Another interpretation of this Lemma consists in the observation
that formulas \re{tP_CT} and \re{GT} are equivalent to the
Christoffel and Geronimus transformations of orthogonal
polynomials which in turn are special cases of the more general
rational spectral (or Darboux) transformations \cite{ZhS},
\cite{BM}.

If the polynomials $P_n(x)$ are orthogonal with respect to a
linear functional $\sigma$: \be \langle \sigma, P_n(x) P_m(x)
\rangle =0, \quad n \ne m , \lab{ort_sigma_P} \ee then the
polynomials $\t P_n(x)$ are orthogonal with respect to the
functional $(x-\theta) \sigma$, i.e. \be \langle \sigma,
(x-\theta) \t P_n(x) \t P_m(x) \rangle =0, \quad n \ne m .
\lab{ort_sigma_tP} \ee

Note also that \be A_n = \frac{P_{n+1}(\theta)}{P_n(\theta)}
\lab{A_PP} \ee which can easily be verified directly from
\re{3-term_AC}.

The next Lemma will be useful in the identification  of polynomial
systems with known families of orthogonal polynomials.

\begin{lem}
Let $P_n(x)$ and $\t P_n(x)$ be two systems of orthogonal
polynomials defined as in the previous Lemma. Assume moreover that
$\theta=\chi^2 \ge 0$.

Define the following monic polynomials \be S_{2n}(x) = P_n(x^2) ,
\quad S_{2n+1}(x) = (x-\chi) \t P_n(x^2). \lab{S_PP} \ee The
polynomials $S_n(x)$ are orthogonal and satisfy the recurrence
relation \be S_{n+1} + (-1)^n \chi S_n(x) + v_n S_{n-1}(x) =
xS_n(x), \lab{rec_S} \ee where \be v_{2n} =-C_n, \; v_{2n+1}=-A_n.
\lab{v_CA} \ee In particular, when $\chi=0$, the polynomials
$S_n(x)$ become symmetric, i.e. $S_{n}(-x) = (-1)^n S_n(x)$.

Conversely, assume that polynomials $S_n(x)$ satisfy the
recurrence relation \re{rec_S} with some coefficients $v_n$ and a
real constant $\chi$. Then the polynomials  $P_n(x)$ and $\t
P_n(x)$ defined by \re{S_PP}, are orthogonal and obey the
recurrence relations \re{3-term_AC} and \re{3-term_tP} with
$\theta=\chi^2$.
\end{lem}
This Lemma is due to Chihara \cite{Chi_Bol}. For further
development and applications of this Lemma see, e.g. \cite{MP},
\cite{VZ_big}.

We shall call the polynomials $\t P_n(x)$, the companion
polynomials with respect to $P_n(x)$ (sometimes the polynomials
$\t P_n(x)$ are referred to as  the kernel polynomials
\cite{Chi}).

The importance of Lemma {\bf 1} lies in the fact that almost all
known classical orthogonal polynomials in the Askey table
\cite{KLS} admit a representation of their recurrence relation in
the form  \re{3-term_AC} with simple explicit coefficients $A_n,
C_n$. The corresponding companion polynomials $\t P_n(x)$ satisfy
similar recurrence relation \re{3-term_tP} and belong to the same
class of orthogonal polynomials (with shifted parameters). For the
polynomials at the top of the Askey hierarchy - the Askey-Wilson
(or $q$-Racah) polynomials,  this was first observed by Chihara in
\cite{Laura}. Further development of this subject in terms of
Darboux transformations can be found in \cite{ZhS_D}.

We see that the Bannai-Ito polynomials admit as well this
representation with the coefficients $A_n, C_n$ given by
\re{A_BIP} and \re{C_BIP}.

We can hence introduce their companion polynomials, which we
denote by $W_n(x)$, as follows \be W_n(x) = \frac{P_{n+1}(x) - A_n
P_n(x)}{x-\rho_1} \lab{W_P} \ee with $A_n$ given by \re{A_BIP}.
The polynomials $W_n(x)$ satisfy the recurrence relation
\re{3-term_tP}. Specifically,  this recurrence relation reads \be
W_{n+1}(x) +(-1)^n \rho_2 W_n(x) + v_n W_n(x) = x W_n(x),
\lab{rec_W} \ee with \ba
&&v_{2n} = -\frac{n(n+\rho_1-r_1+1/2)(n+\rho_1-r_2+1/2)(n-r_1-r_2)}{(2n+1+g)(2n+g)}, \nonumber \\
&&v_{2n+1} =
-\frac{(n+g+1)(n+\rho_1+\rho_2+1)(n+\rho_2-r_1+1/2)(n+\rho_2-r_2+1/2)}{(2n+1+g)(2n+g+2)},
\lab{v_BI} \ea where we denote
$$
g=\rho_1+\rho_2-r_1-r_2.
$$
It is seen that the companion polynomials $W_n(x)$ satisfy a
recurrence relation of the type \re{rec_S} with $\chi=\rho_2$.
Hence a pair of polynomials $U_n(x), \: V_n(x)$ can be constructed
via \be U_n(x^2) =W_{2n}(x), \quad V_n(x^2)=
\frac{W_{2n+1}(x)}{x-\rho_2}, \quad n=0,1,2,\dots . \lab{def_UV}
\ee (It is seen that $W_{2n}(x)$ are polynomials depending only on
$x^2$ and that the polynomials $W_{2n+1}(x)$ all have in common
the factor $x-\rho_2$).

From the above Lemmas, it follows that the polynomials $U_n(y),
V_n(y)$ satisfy the following system of relations \be U_n(y) =
V_n(y) + v_{2n} V_{n-1}(y), \quad (y-\rho_2^2)V_n(y) = U_{n+1}(y)
+ v_{2n+1} U_n(y) \lab{UV_rel} \ee from which the three-term
recurrence relations for the polynomials $U_n(y)$ and $V_n(y)$ can
be derived: \be U_{n+1}(y) +(\rho_2^2+v_{2n} + v_{2n+1}) U_n(y) +
v_{2n} v_{2n-1} U_{n-1}(y) = y U_n(y) \lab{U_rec} \ee and \be
V_{n+1}(y) +(\rho_2^2+v_{2n+2} + v_{2n+1}) V_n(y) + v_{2n}
v_{2n+1} V_{n-1}(y) = y V_n(y). \lab{V_rec} \ee From the explicit
expressions \re{v_BI} of the coefficients $v_n$, it is possible to
conclude that the polynomials $U_n(x^2), \: V_n(x^2)$ are Wilson
polynomials \cite{KLS}. Omitting technical details we have \be
U_n(x^2) = \kappa_n^{(1)} \: {_4}F_3 \left( {-n, n+g+1, \rho_2+x,
\rho_2-x  \atop  \rho_1+\rho_2+   1, \rho_2-r_1+1/2,
\rho_2-r_2+1/2 }   ; 1\right) \lab{U_Wil} \ee and \be V_n(x^2) =
\kappa_n^{(2)} \: {_4}F_3 \left( {-n, n+g+2, \rho_2+1+x,
\rho_2+1-x  \atop  \rho_1+\rho_2+   2, \rho_2-r_1+3/2,
\rho_2-r_2+3/2 }   ; 1\right), \lab{V_Wil} \ee where the
normalization coefficients (making the polynomials monic) are
$$
\kappa^{(1)}_n=
\frac{(1+\rho_1+\rho_2)_n(\rho_2-r_1+1/2)_n(\rho_2-r_1+1/2)_n}{(n+g+1)_n},
\quad \kappa^{(2)}_n=
\frac{(2+\rho_1+\rho_2)_n(\rho_2-r_1+3/2)_n(\rho_2-r_1+3/2)_n}{(n+g+2)_n}
.
$$

(Note that the polynomials $V_n(x^2)$ are obtained from the
polynomials $U_n(x^2)$ by simply shifting the parameter $\rho_2$
by 1: $\rho_2 \to \rho_2+1$.)

We thus have an explicit expression for the polynomials $W_n(x)$
(which are Christoffel transforms of the Bannai-Ito polynomials)
in terms of the Wilson polynomials: \be W_{2n}(x) = U_n(x^2) ,
\quad W_{2n+1}(x) = (x-\rho_2) V_n(x^2), \lab{W_UV} \ee where the
polynomials $U_n(x^2), V_n(x^2)$ are given by \re{U_Wil},
\re{V_Wil}.

Hence, by \re{GT}, we get for the Bannai-Ito polynomials \be
P_n(x) = W_n(x) - C_n W_{n-1}(x), \lab{BI_W} \ee where $C_n$ is
given by \re{C_BIP}. Formula \re{BI_W} can be considered as an
alternative explicit expression of the Bannai-Ito polynomials.
Like in \re{hyp_P_ev} and \re{hyp_P_od}, we have a presentation of
the Bannai-Ito polynomials $P_n(x)$ as a linear combination of two
hypergeometric polynomials ${_4}F_3(1)$. In \re{BI_W} however, the
hypergeometric functions   are -1 - balanced and are thus bona
fide Wilson polynomials. Note that the representation
\re{hyp_P_ev}, \re{hyp_P_od} is equivalent to the one obtained  by
Bannai and Ito for their polynomials by using a direct limit
process from the q-Racah polynomials \cite{BI}. The representation
\re{BI_W} in terms of the Wilson polynomials seems to be new. When
$N \to \infty$, the Bannai-Ito polynomials tend to the big -1
Jacobi polynomials \cite{VZ_big}.  Formula \re{BI_W} then becomes
the formula that gives the expression of the big -1 Jacobi
polynomials in terms of the ordinary Jacobi polynomials
\cite{VZ_big}.

We shall refer to the polynomials $W_n(x)$ as the complementary
Bannai-Ito (CBI) polynomials. They are Christoffel transforms of
the Bannai-Ito polynomials. They have a simpler structure (e.g.
for every $n$ the polynomial $W_n(x)$ is expressed in terms of
only one hypergeometric function ${_4}F_3(1)$).
 However, in contrast to the Bannai-Ito polynomials, the CBI polynomials do not satisfy the Leonard duality property.
 In other words, they are not eigenpolynomials of a Dunkl shift operator of the form
 \re{L_FG}.
 Indeed, the CBI polynomials have recurrence coefficients that do not belong to the BI class, and hence
 they do not arise in the classification of all orthogonal polynomials with the Leonard duality property
 realized by Bannai and Ito \cite{BI}.

\section{Symmetry factor for the Bannai-Ito operator}
\setcounter{equation}{0} Let us return to the Bannai-Ito operator
\re{L_FG}, where $F(x)$ and $G(x)$ are given by \re{FG_simpl}.
Consider the formal Lagrange adjoint $L^*$ of $L$:
$$
L^* = F(x)I - R^* F(x) + (T^+R)^* G(x) -G(x)I.
$$
We will assume that $R^*=R$ and ${T^+}^*=T^-$ (these are natural
and conventional conditions for difference operators on the real
line).

There are in addition, the obvious relations
$$
T^-R=RT^+, \quad  T^+R f(x) = f(-x-1) T^+R,
$$
where $f(x)$ is an arbitrary function of the real variable $x$.
From these relations it follows that the operator $T^+R$ is
formally self-adjoint: $(T^+R)^* = T^+R$.

We then obtain
$$
L^* = F(x) I -F(-x) R +G(-x-1) T^+R -G(x) I.
$$
It is visible that the Bannai-Ito operator is not symmetric, i.e.
$L^* \ne L$ for any choice of the parameters $r_1, r_2, \rho_1,
\rho_2$.

Nevertheless, it is possible to find a {\it symmetry factor}
$\varphi(x)$ for the operator $L$ such that \be (\varphi(x) L)^* =
\varphi(x) L . \lab{sym_cond} \ee An operator L is called
symmetrizable if there exists a real-valued functions $\varphi(x)$
ensuring that property \re{sym_cond} holds. The similar approach
for differential operators and their polynomial solutions is
familiar \cite{LR}.

Let us show that the Bannai-Ito operator is symmetrizable for generic values of the parameters $r_1, r_2, \rho_1, \rho_2$.

We have \be \varphi(x) L = \varphi(x) F(x) - \varphi(x) F(x) R +
\varphi(x) G(x) T^+ R -\varphi(x) G(x) \lab{phi_L} \ee and \be
(\varphi(x) L)^* = \varphi(x) F(x) - \varphi(-x) F(-x) R +
\varphi(-x-1) G(-x-1) T^+ R -\varphi(x) G(x). \lab{phi_L*} \ee
Comparing \re{phi_L*} and \re{phi_L}, we see that condition
\re{sym_cond} is valid iff the following two conditions \be
\varphi(-x)F(-x) = \varphi(x) F(x), \quad \varphi(-x-1)G(-x-1) =
\varphi(x) G(x) \lab{2_cond_FG} \ee are fulfilled.

Clearly, these conditions are equivalent to \be \varphi(x) F(x) =
E_0(x), \quad \varphi(x) G(x) = E_1(x+1/2), \lab{phi_FG_E} \ee
where $E_0(x), E_1(x)$ are even functions, i.e. $E_{0}(-x)=E_0(x),
\; E_{1}(-x)=E_1(x) $. From \re{phi_FG_E} we have \be
\frac{F(x)}{G(x)} = \frac{E_0(x)}{E_1(x+1/2)} \lab{EG_c1} \ee as
well as \be \frac{F(-x)}{G(-x)} = \frac{E_0(x)}{E_1(x-1/2)}
\lab{EG_c2} \ee and \be \frac{F(-x-1)}{G(-x-1)} =
\frac{E_0(x+1)}{E_1(x+1/2)} . \lab{EG_c3} \ee From \re{EG_c1} and
\re{EG_c2}, we get an equation for the unknown function $E_1(x)$
\be \frac{E_1(x+1/2)}{E_1(x-1/2)} = \frac{F(-x) G(x)}{F(x) G(-x)}.
\lab{eq_E1} \ee Similarly, from \re{EG_c1} and \re{EG_c3}, we find
an equation for the unknown function $E_0(x)$: \be
\frac{E_0(x+1)}{E_0(x)}= \frac{F(-x-1) G(x)}{F(x) G(-x-1)} .
\lab{eq_E0} \ee The general solution of equations \re{eq_E0} and
\re{eq_E1} can easily be obtained from the expressions for the
functions $F(x)$ and $G(x)$: \be E_0(x) = \frac{\sigma_0(x)}{x} \:
\frac{\Gamma(1/2-r_1+x)\Gamma(1/2-r_1-x)\Gamma(1+\rho_1+x)\Gamma(1+\rho_1-x)}{\Gamma(1/2+r_2+x)\Gamma(1/2+r_2-x)\Gamma(-\rho_2+x)\Gamma(-\rho_2-x)},
\lab{E_0_sol} \ee \be E_1(x) = \frac{\sigma_1(x)}{x} \:
\frac{\Gamma(1-r_1+x)\Gamma(1-r_1-x)\Gamma(1/2+\rho_1+x)\Gamma(1/2+\rho_1-x)}{\Gamma(r_2+x)\Gamma(r_2-x)\Gamma(1/2-\rho_2+x)\Gamma(1/2-\rho_2-x)},
\lab{E_1_sol} \ee where $\sigma_0(x)$ and $\sigma_1(x)$ are
arbitrary functions with period $1$: $\sigma_{0,1}(x+1) =
\sigma_{0,1}(x)$. The fact that both functions $E_{0}(x)$ and
$E_1(x)$ must be even, requires that both  $\sigma_{0}(x)$ and
$\sigma_1(x)$ be odd:
$$
\sigma_{0,1}(-x) = -\sigma_{0,1}(x) .
$$
From condition \re{EG_c1} it follows that
$\sigma_1(x)=\sigma_0(x)$. Then, from \re{phi_FG_E}, we find the
function $\varphi(x)$ \be \varphi(x) = -2\sigma_0(x)\,{\frac
{\Gamma \left( \rho_{{1}}-x \right) \Gamma \left( x-r_{{1 }}+1/2
\right) \Gamma \left( -x-r_{{1}}+1/2 \right) \Gamma  \left( x+
1+\rho_{{1}} \right) }{\Gamma  \left( x+1-\rho_{{2}} \right)
\Gamma
 \left( x+r_{{2}}+1/2 \right) \Gamma  \left( r_{{2}}+1/2-x \right)
\Gamma  \left( -x-\rho_{{2}} \right) }}. \lab{phi_expr} \ee
Conversely, assume that the function $\varphi(x)$ is given by
\re{phi_expr}. It can then be easily verified that the function
$\varphi(x)$  is a solution of conditions \re{2_cond_FG} iff the
function $\sigma_0(x)$ satisfies the two conditions:

(i) be 1-periodic, i.e. $\sigma_0(x+1)=\sigma_0(x)$;

(ii) be odd, i.e. $\sigma_0(-x) = -\sigma_0(x)$.

Hence we succeeded in constructing the desired symmetry factor
$\varphi(x)$.

It should be noted that the function $\varphi(x)$ can be presented
in several equivalent forms if one exploits the classical formula
$$
\Gamma(z) \Gamma(1-z) =\frac{\pi}{\sin \pi z}.
$$
For example, one can present $\varphi(x)$ in the form \be
\varphi(x) = \sigma_2(x) \: \frac{\Gamma(\rho_1-x)
\Gamma(\rho_1+1+x)\Gamma(\rho_2-x)\Gamma(1+\rho_2+x)
}{\Gamma(r_2+1/2+x) \Gamma(r_2+1/2-x) \Gamma(r_1+1/2+x)
\Gamma(r_1+1/2-x)}, \lab{phi_sym} \ee where the function
$\sigma_2(x)$ has the same properties (i) and (ii) as
$\sigma_0(x)$.

Let us introduce the operator $M=\varphi(x) L$, where $\varphi(x)$
is given by \re{phi_expr}. This operator will thus be formally
self-adjoint, or symmetric:  $M^*=M$.

\section{Bannai-Ito grid and weight function for the Bannai-Ito polynomials}
\setcounter{equation}{0}
One of the most important property of
Dunkl shift operators like $L$ or $M$ is that there is a discrete
set of points on the real line which is invariant under the action
of these operators.

Indeed, let us introduce  the Bannai-Ito (BI) grid of the first
type $x_s$ as a discrete set of real points defined by the two
conditions \be -x_{2s} = x_{2s-1}, \quad   -1-x_{2s} =x_{2s+1},
\quad s=0, \pm 1, \pm 2, \dots. \lab{BI_cond1} \ee It is easy to
verify that the general solution for the grid $x_s$ satisfying
conditions \re{BI_cond1}, is \be x_s = \left\{ a+\frac{s}{2} \quad
\mbox{if $s$} \; \mbox{is even}, \atop -a-\frac{s+1}{2} \quad
\mbox{if $s$} \; \mbox{is odd} , \right . \lab{BI_grid1} \ee where
$a$ is an arbitrary real parameter.

Equivalently, \be x_s = -1/4+(-1)^s (1/4+a+s/2). \lab{BI_grid1_e}
\ee

Using properties \re{BI_cond1}, we can conclude that for any Dunkl
shift operator \be H = A(x) R + B(x) T^+R + C(x) \lab{H_ABC} \ee
with arbitrary functions $A(x), B(x), C(x)$, we have for any
function $f(x)$ \be H f(x_s) = \xi_s f(x_{s+1}) + \eta_s f(x_s) +
\zeta_s f(x_{s-1}), \lab{H-3-diag} \ee where \be \xi_s = \left\{
B(x_s)  \quad \mbox{if $s$} \; \mbox{is even}, \atop      A(x_s)
\quad \mbox{if $s$}  \; \mbox{is odd} ,  \right .  \lab{xi_BIG}
\ee \be \zeta_s = \left\{ A(x_s) \quad \mbox{if $s$} \; \mbox{is
even}, \atop      B(x_s) \quad \mbox{if $s$}  \; \mbox{is odd} ,
\right .   \lab{zeta_BIG} \ee \be \eta_s = C(x_s).  \lab{eta_BIG}
\ee In other words, the Dunkl shift operator becomes a 3-diagonal
matrix (either finite or infinite) in the basis $f(x_s)$.

Similarly, we can define the BI grid of the second type $y_s$ by
the two conditions \be -y_{2s} = y_{2s+1}, \quad   -1-y_{2s}
=y_{2s-1}, \quad s=0, \pm 1, \pm 2, \dots. \lab{BI_cond2} \ee with
the general solution \be y_s = \left\{ b-\frac{s}{2} \quad
\mbox{if $s$} \; \mbox{is even}, \atop -b+\frac{s-1}{2} \quad
\mbox{if $s$} \; \mbox{is odd} , \right . \lab{BI_grid2} \ee where
$b$ is an arbitrary  real parameter.

Equivalently, \be y_s = -1/4+(-1)^s(1/4+b-s/2). \lab{BI_grid2_e}
\ee

We then have that the operator $H$ given by \re{H_ABC} is
3-diagonal as in \re{H-3-diag} with \be \xi_s = \left\{ A(y_s)
\quad \mbox{if $s$} \; \mbox{is even}, \atop      B(y_s) \quad
\mbox{if $s$}  \; \mbox{is odd} ,  \right .  \lab{xi_BIG2} \ee \be
\zeta_s = \left\{ B(y_s) \quad \mbox{if $s$} \; \mbox{is even},
\atop A(y_s) \quad \mbox{if $s$}  \; \mbox{is odd} , \right .
\lab{zeta_BIG2} \ee \be \eta_s = C(y_s).  \lab{eta_BIG2} \ee

In particular, the BI difference equation \re{LPP} for the BI
polynomials can be presented in the form \be \xi_s P_n(x_{s+1}) +
\eta_s P_n(x_s) + \zeta_s P_n(x_{s-1}) = \lambda_n P_n(x_s),
\lab{BI_3DFR} \ee where the coefficients $\xi_s, \eta_s, \zeta_s$
are given by \re{xi_BIG}-\re{eta_BIG} (or by
\re{xi_BIG2}-\re{eta_BIG2}) with
$$A(x)=-F(x), \; B(x)=G(x), \; C(x)=F(x)-G(x) .$$
In turn, property \re{BI_3DFR} means that the BI polynomials
possess the Leonard duality property \cite{BI}, \cite{Ter}: they
satisfy both a 3-term recurrence relation with respect to $n$ and
a 3-term difference equation with respect to the argument index
$s$. Bannai and Ito found all orthogonal polynomials satisfying
the Leonard duality property \cite{BI}. The corresponding
classification, referred to as the  Bannai-Ito theorem, is a
generalization of the Leonard theorem \cite{Leonard} that deals
only with polynomials orthogonal on a finite set of points. We
have thus derived the Leonard duality property of the BI
polynomials starting from the Dunkl shift difference equation
\re{LPP}.

Note that the Bannai-Ito grid satisfies the linear equation \be
x_{s+1} + x_{s-1} + 2 x_s +1 =0 \lab{BI_grid_eq} \ee which is a
special case of the generic equation \be x_{s+1} + x_{s-1} -
(q+q^{-1}) x_s  = const \lab{AW_grid_eq} \ee characterizing the
Askey-Wilson grids \cite{Ismail}, \cite{NSU}. As expected, the
Bannai-Ito grid corresponds to the case $q=-1$ of the Askey-Wilson
grid.

Consider now the special case when the operator $H$ is symmetric
$H^*=H$, with real functions $A(x),B(x),C(x)$. This is equivalent
to the conditions \be A(-x)=A(x), \quad  B(-x-1)=B(x) .
\lab{AB_sym} \ee It is easily verified that in this case the
3-diagonal matrix in \re{H-3-diag} becomes symmetric as well, i.e.
$\xi_s=\zeta_{s+1}$, or in details: \be H f(x_s) = \zeta_{s+1}
f(x_{s+1}) + \eta_s f(x_s) + \zeta_s f(x_{s-1}).
\lab{H-3-diag_sym} \ee This result does not depend on the type of
the BI grid. We already know that the operator $M=\varphi(x) L$ is
symmetric $M^* = M$. Assume further that the function $\varphi(x)$
is positive, i.e $\varphi(x)>0$ for some interval $x_0<x<x_1$. We
can then introduce the operator \be \t L =\varphi^{1/2}(x) L
\varphi^{-1/2}(x). \lab{def_tL} \ee By condition \re{sym_cond},
this operator will be symmetric $\t L^* = \t L$ on the interval
$[x_0,x_1].$

Assume additionally that the corresponding 3-diagonal matrix is
finite, say it has dimension $(N+1) \times (N+1)$. This is
equivalent to the conditions \be \zeta_0=\zeta_{N+1}=0
\lab{zeta_0} \ee and $\zeta_s>0, \; s=1,2,\dots,N$. Let $q_n(s)$
be a set of eigenvectors of the s-diagonal matrix $\t L$: \be
\zeta_{s+1} q_n(s+1) + \eta_s q_n(s) + \zeta_s q_n(s) = \lambda_n
q_n(s), \quad s,n=0,1,2,\dots,N . \lab{eigen_q} \ee The operator
$\t L$ is a Hermitian matrix in the basis $q_n(s)$. By elementary
linear algebra, it is well known that the coordinates $q_n(s)$
(under appropriate normalization) satisfy the orthonormality
property \be \sum_{s=0}^N q_n(s)q_m(s) = \delta_{nm} . \lab{ort_q}
\ee From the definition \re{def_tL} of the operator $\t L$, it
then follows that \be q_n(s) = \kappa^{-1}_n
\varphi^{1/2}(x_s)P_n(x_s) \lab{qP} \ee with some normalization
factor $\kappa_n$.

We thus arrive on the one hand, at the orthogonality relation for
the BI polynomials: \be \sum_{s=0}^N w_s P_n(x_s) P_m(x_s) =
\kappa_n^2 \: \delta_{nm}, \lab{ort_wPP} \ee where the discrete
weights are defined by the formula \be w_s = \varphi(x_s) .
\lab{w_s_phi} \ee By construction, the weights are positive
$w_s>0$ as should be for positive definite orthogonal polynomials
$P_n(x)$.

On the other hand, we have formula \re{N-ort} giving the
orthogonality relation of the BI polynomials in the finite case,
where $x_s$ are simple roots of the polynomial $P_{N+1}(x)
=(x-x_0)(x-x_1) \dots(x-x_N)$. This means that the roots $x_s$ can
be parametrized in terms of the BI grid.

Explicitly, the roots $x_s$ of the polynomial $P_{N+1}(x)$ can be
found from formula \re{BI_W} which expresses the BI polynomials in
terms of the Wilson polynomials. Omitting technical details, we
present the results.

When $N$ is even and condition \re{ev_N_rr} is assumed then \be
x_s = -1/4+(-1)^s(\rho_2+1/4 +s/2), \quad s=0,1,2,\dots, N,
\lab{x_s_even_N} \ee i.e. the roots for the case of even $N$
belong to the BI grid of the first type \re{BI_grid1_e} with
$a=\rho_2$.

When $N$ is odd and condition $r_1 + r_2 =(N+1)/2$ is assumed then
\be x_s = -1/4+(-1)^s(-1/4+r_1-s/2), \quad s=0,1,2,\dots, N,
\lab{x_s_odd_N} \ee i.e. the roots for the case of odd $N$ belong
to the BI grid of the second type \re{BI_grid2_e} with
$b=r_1-1/2$.

So for even $N=2,4,6,\dots$ the orthogonality relation reads as in
\re{ort_wPP} with the spectral points $x_s$ given by
\re{x_s_even_N} and the discrete weights expressed as in
\re{w_s_phi}, and where $\varphi(x)$ can be taken in the form
\re{phi_expr}.

For odd $N=1,3,5,\dots$ the spectral points $x_s$ are given by
\re{x_s_even_N} and the discrete weights expressed as in
\re{w_s_phi}, where $\varphi(x)$ can now be taken of the form
\re{phi_sym} (this form is preferable for the case of odd $N$ in
order to avoid singularities when $x=x_s$).

It is easily verified that in all cases the function
$\sigma_0(x_s)$ and $\sigma_2(x_s)$ are constant sign-alternating
sequence, e.g.
$$
\sigma_0(x_s) =(-1)^s const, \quad s=0,1,2,\dots,N
$$

Note finally that the set of spectral points $x_s$ is the union of
two discrete subsets $x_{2s}$ and $x_{2s+1}$ corresponding to even
and odd values of $s$. These two subsets never overlap if the
positivity conditions form the BI polynomials are fulfilled.

\section{Bannai-Ito and complementary Bannai-Ito polynomials as limits $q \to -1$ of the Askey-Wilson polynomials}
\setcounter{equation}{0} The Bannai-Ito polynomials were
introduced in \cite{BI} through a direct limit process from the
$q$-Racah polynomials as $q \to -1$.

We show in this section how the Bannai-Ito polynomials can be
obtained from the Askey-Wilson (AW) polynomials as $q \to -1$. We
also indicate how the complementary Bannai-Ito polynomials can
equally be derived from the AW-polynomials in the same limit.

Consider the Askey-Wilson polynomials $R_n(x(z);a,b,c,d)$ \cite{AW}, \cite{KLS}
\be
R_n(z;a,b,c,d)= \kappa_n^{(1)} \; {_4}\Phi_3 {\left( {q^{-n}, abcd q^{n-1}, az, az^{-1}   \atop ab, ac, ad   } \left | \right . q;q \right)}
\lab{AWP} \ee
where
$$
\kappa_n^{(1)} = a^{-n} (ab,ac,ad;q)_n .
$$
The polynomials $R_n(x(z);a,b,c,d)$ depend on the argument
$x=(z+z^{-1})/2$ and on 4 complex parameters $a,b,c,d$. They
satisfy the recurrence relation \cite{KLS} \be A_n R_{n+1}(x) +
(a+a^{-1}-A_n-C_n)R_n(x) +C_nR_{n-1}(x) = (z+z^{-1})R_n(x)= 2x
R_n(x) \lab{REC_AW} \ee with coefficients \ba
&&A_n = \frac{(1-abq^{n})(1-acq^{n})(1-adq^{n})(1-abcdq^{n-1})}{a(1-abcd q^{2n-1})(1-abcd q^{2n})}, \nonumber \\
&&C_n =
\frac{a(1-q^{n})(1-bcq^{n-1})(1-bdq^{n-1})(1-cdq^{n-1})}{(1-abcd
q^{2n-1})(1-abcd q^{2n-2})} . \lab{AC_AW} \ea In order to recover
the BI-polynomials we choose the following parametrization \be
q=-e^{\ve}, \; a=-ie^{\ve \alpha}, \; b=-ie^{\ve \beta}, \;
c=ie^{\ve \gamma}, \; d=ie^{\ve \delta}, \lab{AW_BI_param} \ee
where $\alpha, \beta, \gamma, \delta$ are real parameters. In what
follows these parameters are more conveniently expressed in terms
of the 4 parameters $r_{1,2}, \rho_{1,2}$: \be \alpha=2 \rho_1
+1/2, \quad \beta=2 \rho_2 +1/2, \quad \gamma=-2 r_2 +1/2, \quad
\delta=-2 r_1 +1/2 . \lab{al_bet_rr} \ee Moreover, we take \be
z=ie^{-2\ve y} \lab{z_y} \ee with $y$ a real parameter.

It can be checked that the limit $q \to -1$ of the Askey-Wilson
polynomials \be \lim_{q \to -1} R_n(x(z)) = R^{(-1)}(y)
\lab{limit_RR} \ee does exist, and that $R_n^{(-1)}(y)$ is a
polynomial of degree $n$ in the argument $y$.

Dividing the recurrence relation \re{REC_AW} by $1+q$ and taking
afterwards the limit $\ve \to 0$ (which is equivalent to the limit
$q \to -1$), we obtain a recurrence relation for the limiting
polynomials $R_n^{(-1)}(y)$: \be A^{(-1)}_n R_{n+1}^{(-1)}(y)
+(1/4+\rho_1 -A^{(-1)}_n - C^{(-1)}_n ) R_n^{(-1)}(y) + C_n^{(-1)}
R_{n-1}^{(-1)}(y) = y R_n^{(-1)}(y), \lab{rec_R-1} \ee where the
coefficients $A^{(-1)}_n, \: C^{(-1)}_n$ have expressions
coinciding with \re{A_BIP} and \re{C_BIP}. The polynomials
$R_n^{(-1)}(y)$ are thus identified with the Bannai-Ito
polynomials (up to a shift of the argument): \be \hat
R_n^{(-1)}(y) = P_n(y-1/4;r_1,r_2,\rho_1,\rho_2). \lab{R_BIP} \ee
$\hat R_n^{(-1)}(y)$ in \re{R_BIP} denotes the monic version of
the polynomials $R_n^{(-1)}(y)$, while
$P_n(y;r_1,r_2,\rho_1,\rho_2)$ stands for the monic Bannai-Ito
polynomials defined in \re{rec_AC_BIP}. The parameters
$r_1,r_2,\rho_1, \rho_2$ for the Bannai-Ito polynomials stem under
the limiting procedure \re{al_bet_rr}, from the corresponding 4
parameters of the Askey-Wilson polynomials.

For the complementary Bannai-polynomials, we use the same
procedure but with a slightly different parametrization: \be
q=-e^{\ve}, \; a=ie^{\ve (2 \rho_1 +3/2)}, \; b=-ie^{\ve (\rho_2
+1/2)}, \; c=ie^{\ve (-2 r_2 +1/2)}, \; d=ie^{\ve (-2 r_1
+1/2)},\; z=ie^{-2\epsilon y} . \lab{AW_preBI_param} \ee The main
difference between \re{AW_BI_param} and \re{AW_preBI_param} is a
change of sign in $a$ and the shift $\rho_1 \to \rho_1+1/2$.

It is directly verified that in the limit $\ve \to 0$, we get the
recurrence relation \re{rec_W} that defines the  CBI polynomials
and therefore that  \be \lim_{\ve \to 0} R_n(z;a,b,c,d) =
W_n(y-1/4;r_1,r_2,\rho_1,\rho_2). \lab{lim_AW_BI} \ee Thus the BI
and CBI polynomials are obtained through very similar $q \to -1$
limits of the Askey-Wilson polynomials.

\section{The Bannai-Ito Dunkl shift  operator as a limiting form of the Askey-Wilson difference operator}
\setcounter{equation}{0} In this section, we describe how the
Bannai-Ito Dunkl shift operator appears in the limiting process $q
\to -1$ from the difference Askey-Wilson operator.

The Askey-Wilson polynomials $R_n(x(z))$ satisfy the difference
equation \cite{KLS} \be \Omega(z) R_n(x(zq)) + \Omega(z^{-1})
R_n(x(zq^{-1})) -(\Omega(z) + \Omega(z^{-1})) R_n(x(z)) =
\Lambda_n R_n(x(z)), \quad n=0,1,2,\dots , \lab{AW_DE}  \ee where
$$
\Lambda_n = (q^{-n}-1)(1-abcd q^{n-1})
$$
is the eigenvalue and \be \Omega(z) =
\frac{(1-az)(1-bz)(1-cz)(1-dz)}{(1-z^2)(1-qz^2)} . \lab{Om_AW} \ee
The linear operator in the lhs of \re{AW_DE} (acting on the
variable $z$) is the Askey-Wilson difference operator \cite{AW},
\cite{KLS}.

Consider the limiting form of the difference equation \re{AW_DE}
when $q \to -1$. As in the previous section, we choose the
parametrization \re{al_bet_rr} and \re{z_y}. We already showed
that the Askey-Wilson polynomials $R_n(x(z))$ become the
Bannai-Ito polynomials $R_n^{(-1)}(y)$. The operation $R_n(z) \to
R_n(zq)$ is reduced to the operation $R_n^{(-1)}(y) \to
R_n^{(-1)}(-y+1/2)$ while the operation $R_n(z) \to R_n(zq^{-1})$
becomes the operation $R_n^{(-1)}(y) \to R_n^{(-1)}(-y-1/2)$.

It is then straightforward to verify that in the limit $\ve \to 0$
we get the equation \be \Phi_1(y) R_n^{(-1)}(-y+1/2) + \Phi_2(y)
R_n^{(-1)}(-y-1/2) -(\Phi_1(y) + \Phi_2(y)) R_n^{(-1)}(y) =
\lambda_n R_n^{(-1)}(y), \lab{Phi_R} \ee where \ba &&\Phi_1(y) =
\lim_{\ve \to 0} \frac{\Omega(z)}{4(1+q)} =
-\frac{2(y-\rho_1-1/4)(y-\rho_2-1/4)}{4y-1}, \nonumber \\
&&\Phi_2(y) = \lim_{\ve \to 0} \frac{\Omega(z^{-1})}{4(1+q)} =
\frac{2(y-r_1+1/4)(y-r_2+1/4)}{4y+1} \lab{Phi_12} \ea and that the
eigenvalue \be \lambda_n = \lim_{\ve \to 0}
\frac{\Lambda_n}{4(1+q)} \lab{lam_lim} \ee has the expression
\re{lambda_fix}.

This gives the Bannai-Ito polynomials  $R_n^{(-1)}(y)$  as
eigenfunctions of some Dunkl shift operator. In operator form this
eigenvalue equation can be presented as follows:
 \be H R_n^{(-1)}(y) = \lambda_n
R_n^{(-1)}(y). \lab{HRR} \ee $H$ in \re{HRR}, acts on the space of
functions $f(y)$ of argument $y$ according to \be H = \Phi_1(y)
T^{-1/2}R + \Phi_2(y) T^{1/2}R -(\Phi_1(y) + \Phi_2(y)),
\lab{H_Dunkl} \ee and  $T^{h}$ stands for the shift operator
$$
T^{h} f(y) = f(y+h).
$$
The Dunkl shift  operator $H$ contains (apart from the identity
operator) the operators $T^{1/2}, T^{-1/2}$ and $R$.

We can now reduce the operator \re{H_Dunkl} to the "canonical"
form which only involves the operators $T^+$ and $R$. It goes as
follows. The similarity transformation
$$
\t H = T^{h} H T^{-h}
$$
with an arbitrary real $h$ is a unitary transformation of the
operator $H$ (recall that by $T^{+}$ and $T^{-}$ we mean the
operators $T^{+1}$ and $T^{-1}$ respectively). Under this
transformation $H$ goes into \be \t H = \Phi_1(y+h)T^{2h-1/2}R +
\Phi_2(y+h)T^{2h+1/2}R - \Phi(y+h) - \Phi_1(y+h). \lab{tH_can} \ee
For $h=1/4$, we have \be \t H = \Phi_1(y+1/4)R +
\Phi_2(y+1/4)T^{+}R - \Phi_1(y+1/4) - \Phi_2(y+1/4),
\lab{H_shift1} \ee where \be \Phi_1(y+1/4) =
-\frac{(y-\rho_1)(y-\rho_2)}{2y}=-F(y), \quad \Phi_2(y+1/4) =
-\frac{(y-r_1+1/2)(y-r_2+1/2)}{2y+1}=G(y). \lab{H_shift2} \ee The
functions $G(x), F(x)$ in \re{H_shift2} are precisely those of
\re{FG_simpl}. We have thus reduced the operator $H$ to our
initial "canonical"  Dunkl shift operator $L$ defined in
\re{L_FG}.

Note that there is one more "canonical" form of the Dunkl shift
operator where only $T^-$ and $R$ appear. This other form is
obtained by choosing $h=-1/4$, so that  \be \t H=
\Phi_1(y-1/4)T^{-} R + \Phi_2(y-1/4) R - \Phi_1(y-1/4) -
\Phi_2(y-1/4). \lab{-canon} \ee We see that the  Dunkl shift
operator can be presented in 3 "canonical" form: the "$+$ form"
given by \re{tH_can}, the "$-$ form" given by \re{-canon} and the
"symmetric" form given by \re{H_Dunkl}. Clearly, all these forms
are equivalent and correspond to shifts $y \to y+h$ in the
argument of the polynomials.

As was already mentioned, the complementary BI polynomials do not
satisfy a difference equation of Dunkl shift type. Indeed, under
the choice \re{AW_preBI_param}, the corresponding Askey-Wilson
difference operator does not survive in the limit $q \to -1$.

\section{Symmetric BI polynomials and continuous dual Hahn polynomials}
\setcounter{equation}{0}
We already saw that for generic BI polynomials with real parameters $r_1, r_2, \rho_1, \rho_2$
it is impossible to obtain positive definite polynomials for all $n=0,1,2,\dots$.

There is however, a special case for which it is possible to
construct such positive definite polynomials after the following
change of argument: $x \to ix$. This case arises for symmetric BI
polynomials.

Recall that polynomials $S_n(x)$ are said to be the symmetric if
they satisfy the condition \be S_n(-x) =(-1)^n \: S_n(x) .
\lab{symm_P} \ee Orthogonal polynomials $P_n(x)$ obeying a
recurrence relation of the form
$$
P_{n+1}(x) + b_n P_n(x) + u_n P_{n-1}(x) = xP_n(x)
$$
are symmetric iff $b_n=0$ \cite{Chi} (i.e. iff the diagonal
recurrence coefficient $b_n$ vanishes for all $n=0,1,2,\dots$).

If $b_n=b =const,\; n=0,1,2,\dots$ it is then the "shifted"
polynomials $\t P_n(x) = P_n(x+b)$ that are symmetric.

Let us determine when the (shifted) BI polynomials $R_n^{(-1)}(x)$
(defined by \re{R_BIP}) are symmetric. This is equivalent to
determining when the condition $b_n=-1/4, \; n=0,1,2,\dots$ for
the diagonal recurrence coefficient of the BI polynomials, is
satisfied. It is easy to see that this condition holds iff

either

(i) $\rho_1 = -r_1, \; \rho_2 = -r_2$;

or

(ii) $\rho_1=-r_2, \; \rho_2=-r_1.$

In both cases we have \be b_n = -1/4, \quad n=1,2,3,\dots
\lab{b_sym} \ee and \be u_n=\left\{ -\frac{n(n-4(r_1+r_2))}{16} ,
\quad \mbox{if} \quad n \quad         \mbox{even} \atop
-\frac{(n-4r_1)(n-4r_2)}{16} , \quad \mbox{if} \quad n \quad
\mbox{odd} \right . . \lab{uBI_sym} \ee The corresponding
symmetric polynomials $S_n(x) = \hat R_n^{(-1)}(x)=P_n(x-1/4)$ are
expressed in terms of dual Hahn polynomials. In order to obtain
positive definite polynomials, we need to perform the change of
variable $x \to ix$. This leads to symmetric polynomials with
positive coefficients $u_n$. Indeed, for the new polynomials $\t
S_n(x) = i^{-n}S_n(ix)$, we have the 3-term recurrence relation
\be \t S_{n+1}(x) + \t u_n \t S_{n-1}(x) = x \t S_n(x)
\lab{rec_S_H} \ee with coefficients $\t u_n = i^2 u_n=-u_n$ and
thus \be \t u_n=\left\{ \frac{n(n-4(r_1+r_2))}{16} , \quad
\mbox{if} \quad n \quad \mbox{even} \atop
\frac{(n-4r_1)(n-4r_2)}{16} , \quad \mbox{if} \quad n \quad
\mbox{odd} \right . . \lab{uBI_sym_pos} \ee It is seen that under
the restrictions
$$
r_1<1/4, \quad r_2 <1/4,
$$
one has $\t u_n>0$ for all $n=1,2,3,\dots$. This means that the
polynomials $\t S_n(x)$ are orthogonal on the real axis with a
positive measure.

Comparing the recurrence coefficients \re{uBI_sym_pos} with those
for the continuous dual Hahn polynomials \cite{KLS} and using
Chihara's standard procedure to relate symmetric and non-symmetric
polynomials \cite{Chi}, we can derive an explicit expression for
the polynomials $\t S_n(x)$: \be \t S_{2n}(x) = \kappa_n^{(0)} \:
{_3} F_2 \left({-n, a+2ix, a-2ix   \atop a, a+b} ;1 \right) ,
\quad \t S_{2n+1}(x) = \kappa_n^{(1)} \: x \: {_3} F_2 \left({-n,
a+2ix, a-2ix   \atop a+1, a+b}   ;1 \right), n=0,1,2,\dots
,\lab{S_Hahn} \ee where
$$
a=-2 r_1 + 1/2, \; b=-2 r_2 + 1/2
$$
and the normalization coefficients are
$$
\kappa_n^{(0)}=(-1)^n 2^{-2n} \: (a)_n \: (a+b)_n, \quad \kappa_n^{(1)}=(-1)^n 2^{-2n} \: (a+1)_n \: (a+b)_n .
$$
Recall that the (unnormalized) continuous dual Hahn polynomials
$C_n(x^2;a,b,c)$ depend on 3 parameters $a,b,c$ and are expressed
as \cite{KLS} \be C_n(x^2;a,b,c) = {_3}F_2\left(  {-n, a+ix, a-ix
\atop  a+b,a+c}; 1 \right) . \lab{dual_Hahn} \ee We thus see that
$$
\t S_{2n}(x) = \kappa_n^{(0)} \: C_n(4 x^2;a,b,0), \quad \t
S_{2n+1}(x) = \kappa_n^{(1)} \: x C_n(4 x^2;a,b,1).
$$
From the expression of the weight function of the continuous dual
Hahn polynomials \cite{KLS}, we reconstruct the weight function of
the symmetric BI polynomials \be w(x) = \left |
\frac{\Gamma(a+2ix) \Gamma(b+2ix)}{\Gamma(1/2+2ix)} \right |^2,
\lab{w_dual_H} \ee for which they are orthogonal on the whole real
line
$$
\int_{-\infty}^{\infty} w(x) \t S_n(x) \t S_m(x) dx = 0, \quad  n
\ne m .
$$
The difference equation for these polynomials follows immediately
from the corresponding difference equation  \re{Phi_R} for the BI
polynomials after replacing $y \to iy$: \be \Phi_1(y) \t S_n(-y
-i/2) + \Phi_2(y) \t S_n(-y +i/2) -(\Phi_1(y) +\Phi_2(y)) \t
S_n(y) = \lambda_n \: \t S_n(y), \lab{dfrc_S} \ee where
$$
\Phi_1(y) = -i\frac{(y+ia/2)(y+ib/2)}{2(y +i/4)}, \quad \Phi_2(y) = \Phi_1^*(y)= i\frac{(y-ia/2)(y-ib/2)}{2(y -i/4)}
$$
and the eigenvalue is \be \lambda_n = \left\{ \frac{n}{2} \quad
\mbox{if $n$}  \; \mbox{is even} \atop 2(r_1+r_2) -
\frac{(n+1)}{2} \quad \mbox{if $n$}  \; \mbox{is odd} .  \right .
\lab{lambda_S} \ee Consider the special case when one of the
parameters $r_1$ or $r_2$ vanishes, say $r_2=0$. The recurrence
coefficients $\t u_n$ then become \be \t u_n = \frac{n(n-4
r_1)}{16}, \quad n=1,2,3,\dots \lab{u_MP} \ee These correspond to
the recurrence coefficients of the symmetric Meixner-Pollaczek
polynomials \cite{KLS} \be \t S_n(x) = i^n 2^{-2n} \: (2a)_n \:
{_2}F_1 \left(    {-n, a+2ix  \atop  2a }; 2 \right). \lab{MP_S}
\ee The corresponding Dunkl shift eigenvalue equation is \be
-\frac{i}{4}(2y+ia) \t S_n(-y-i/2) + \frac{i}{4}(2y-ia) \t
S_n(-y+i/2) + \frac{a}{2} \: \t S_n(y) = \lambda_n \: \t S_n(y).
\lab{Dunkl_MP} \ee Taking into account  the property $\t S_n(-y)
=(-1)^n S_n(y)$ of the symmetric polynomials, equation
\re{Dunkl_MP} can be cast in the form \be -\frac{i}{4}(2y+ia) \t
S_n(y+i/2) + \frac{i}{4}(2y-ia) \t S_n(y-i/2) =\frac{a+n}{2} \t
S_n(y) \lab{ord_PM_eq} \ee which coincides with the ordinary
difference equation of the Meixner-Pollaczek polynomials on a
uniform imaginary grid \cite{KLS}.

The weight function in this case is obtained from the weight
function \re{w_dual_H} by putting $b=1/2$:
$$
w(x)=\left | \Gamma(a+2ix)\right |^2,
$$
which coincides  with the well known expression for the weight
function of the symmetric Meixner-Pollaczek polynomials
\cite{KLS}.

The symmetric Meixner-Pollaczek polynomials therefore belong, as a
special case, to the class of BI polynomials. The eigenvalue
equation for these polynomials can be presented either in the form
\re{Dunkl_MP} as an eigenvalue equation for a Dunkl shift
difference operator of first order, or in the form \re{ord_PM_eq}
which is a standard  Sturm-Liouville difference equation of the
second order.

\section{The big and little -1 Jacobi polynomials as limit cases}
\setcounter{equation}{0} Changing the variable to $x=y/h$ with $h$
a real parameter, we can rewrite the BI operator  in the form \be
L f(y) = \frac{(y-\rho_1 h)(y-\rho_2 h)}{2h y}(f(y)-f(-y)) +
\frac{(y-r_1h+h/2)(y-r_2h +h/2)}{h(2y+h)}(f(-y-h)-f(y)). \lab{L_y}
\ee Choose now the following parametrization \be r_1 = a_1/h, \;
r_2 = a_2/h, \; \rho_1 = a_1/h+ b_1, \; \rho_2 = a_2/h+b_2,
\lab{par_r_ab} \ee where the real parameters $a_i,b_i, i=1,2$ do
not depend on $h$. Taking the limit $h \to 0$, we then arrive at
the operator \be L_0=\lim_{h \to 0} L= \frac{(y-a_1)(y-a_2)}{2y}
\partial_y R - \frac{F(y)}{4 y^2}(I-R), \lab{L_dR} \ee where
$$
F(y) = (2b_1+2b_2+1)y^2 -2(a_1 b_2 + a_2b_1)y -a_1a_2 .
$$
The operator \re{L_dR} is a first order  differential operator of
Dunkl type. It was shown in \cite{VZ_Bochner} that the operator
$L_0$ as given in \re{L_dR} with $a_1, a_2,b_1$ and $b_2$ real,
represents the most general first order differential operator of
Dunkl type that has orthogonal polynomials as eigenfunctions. In
general the operator $L_0$ has 4 free parameters. However, if $a_2
\ne 0$, by scaling the argument $y \to \gamma y$ one can set
$a_2=1$ and only 3 independent real parameters $a_1,b_1,b_2$
remain.

Using a slightly different parametrization, we thus have
\cite{VZ_big} \be L_0=g_0(y)(R-I) + g_1(y) \partial_y R,
\lab{L_0_big} \ee where \be g_0(y)= \frac{(\alpha+\beta+1)y^2
+(c\alpha -\beta)y + c}{y^{2}}, \quad
g_1(y)=\frac{2(y-1)(y+c)}{y}. \lab{g_01} \ee The 3 independent
real parameters of the operator $L_0$ are now $\alpha,\beta$ and
$c$. The polynomial eigensolutions $P_n(y), n=0,1,2,\dots$ of
$L_0$ are expressed in terms of the big -1 Jacobi polynomials (see
\cite{VZ_big} for details).

If $a_2=0$, the operator $L_0$ can be reduced to \be
L_0=2(1-y)\partial_y R +(\alpha+\beta+1-\alpha y^{-1})(1-R)
\lab{little_L} \ee with only two independent real parameters
$\alpha, \beta$. The polynomial eigensolutions $P_n(y),
n=0,1,2,\dots$ of the operator \re{little_L} coincide with the
little -1 Jacobi polynomials \cite{VZ_little}.

In \cite{VZ_big} it was shown how the recurrence coefficients of
the big -1 Jacobi polynomials could be recovered through a special
large $N$ limit, from those of the Bannai-Ito polynomials.
Clearly, the eigenvalue problem for the big -1 Jacobi polynomials
can also be obtained from the corresponding eigenvalue problem for
the BI polynomials. For $q \to 1$ this corresponds to the
transition from the q-Racah to the big q-Jacobi polynomials
\cite{K_Racah}.

\section{Conclusion}
\setcounter{equation}{0} The study of the polynomial
eigenfunctions of first order Dunkl shift operator has thus led
synthetically to a thorough characterization of the Bannai-Ito
polynomials: weight function, structure and recurrence
coefficients etc. It has brought to light a relation between
certain Jordan algebras, their representations and orthogonal
polynomials. It made natural the introduction of the complementary
Bannai-Ito polynomials, their expression in terms of Wilson
polynomials and their connection to the BI polynomials. It further
allowed to see that the symmetric BI polynomials are related to
the  dual Hahn polynomials and to study finally, various limiting
cases.

Note that polynomial eigenfunctions of certain linear operators
with reflections were previously considered in the literature (see
for instance \cite{KB}); these solutions however exhibit
non-standard orthogonality properties. In contrast, the polynomial
eigensolutions that have arisen in our treatment, are clearly
orthogonal in the ordinary sense.

\bigskip\bigskip
{\Large\bf Acknowledgments}

\bigskip

The authors are indebted to a referee for useful remarks.
The authors would like to gratefully acknowledge the hospitality
extended to LV and AZ by Kyoto University and to ST and LV by the
Donetsk Institute for Physics and Technology  in the course of
this investigation. The research of LV is supported in part by a
research grant from the Natural Sciences and Engineering Research
Council (NSERC) of Canada.

%\noindent AZ thanks CRM (U de Montr\'eal) for its hospitality.

\newpage

\bb{99}

\bi{AK} M.~Arik and  U.~Kayserilioglu,{\it The anticommutator spin
algebra, its representations and quantum group invariance}, Int.
J. Mod. Phys. A, {\bf 18} (2003), 5039--5046;
arXiv:hep-th/0212127v1.

%\bi{biask} R. Askey, Comments in: "Gabor Szeg\H{o}:Collected Papers",
%Birkh\"auser, Basel, 1982, v.1, 303-305.

%\bi{AI} R. A. Askey and M. E. H. Ismail, {\it Recurrence relations, continued fractions and orthogonal polynomials}, Mem. AMS, {\bf 49} (1984), No. 300,  1--108.

\bi{AW} R.~Askey and J.~Wilson, {\it Some basic hypergeometric orthogonal
polynomials that generalize Jacobi polynomials}, Mem. Amer. Math. Soc. {\bf
54}, No. 319, (1985), 1-55.

%\bi{Atia} M. J. Atia, {\it An example of non-symmetric
%semi-classical form of class s=1. Generalization of a case of
%Jacobi sequence}. Int. J. Math. and Math. Sci. {\bf 24} (10)
%(2000), 673--689.

%\bi{Atia2} M.J.Atia, {\it Some generalized Jacobi polynomials},
%OPSFA 2009, http://wis.kuleuven.be/OPSFA/bookletOPSFA.pdf

%\bi{Bax} G. Baxter, {\it Polynomials defined by a difference
%system}, J.Math.Anal.Appl. {\bf 2} (1961), 223-263.

\bi{BI} E.~Bannai and T.~Ito, {\it Algebraic Combinatorics I:
Association Schemes}. 1984. Benjamin \& Cummings, Mento Park, CA.

\bi{BM} M.I.~Bueno and  F.~Marcell\'an, {Darboux transformation
and perturbation of linear functionals}, Lin. Alg. Appl. {\bf 384}
(2004), 215-�242.

%\bi{Cheikh} Y. Ben Cheikh, {\it Characterization of the Dunkl-classical symmetric orthogonal polynomials}, Appl. Math. and Comput.
%{\bf 187}, (2007) 105--114.

%\bi{Belm} S. Belmehdi, {\it Generalized Gegenbauer polynomials}, J. Comput. Appl Math. {\bf 133} (2001), 195-–205.

%\bi{Cheikh} Y. Ben Cheikh and M.Gaied, {\it Characterization of the Dunkl-classical symmetric orthogonal polynomials}, Appl. Math. and Comput.
%{\bf 187}, (2007) 105--114.

%\bi{Chi2} T.Chihara, {\it Orthogonal polynomials with Brenke type
%generating functions}, Duke Math. J. {\bf 35} (1968), 505--517.

\bi{Chi_Bol} T.~Chihara, {\it On kernel polynomials and related
systems}, Boll. Unione Mat. Ital., 3 Ser., {\bf 19} (1964),
451--459.

\bi{Chi} T.~Chihara, {\it An Introduction to Orthogonal
Polynomials}, Gordon and Breach, NY, 1978.

\bi{Laura} Laura M.~Chihara, {\it Askey-Wilson Polynomials, Kernel
Polynomials and Association Schemes}, Graphs and Combinatorics,
{\bf 9} (1993), 213--223.

%\bi{Chouchene} F.Chouchene,{\it Harmonic analysis associated with the Jacobi-Dunkl operator on $[-\frac{\pi}{2}, \frac{\pi}{2}]$}, J. %Computat. Appl. Math.
%{\bf 178}, (2005), 75-89.

%\bi{Chi_Chi} L.M.Chihara,T.S.Chihara, {\it A class of nonsymmetric
%orthogonal polynomials}. J. Math. Anal. Appl. {\bf 126} (1987),
%275--291.

\bi{Curtin1} B.~Curtin, {\it Spin Leonard pairs}, Ramanujan J.,
{\bf 13} (2007), 319-�332

\bi{Curtin2} B.~Curtin, {\it Modular Leonard triples}, Linear
Algebra and its Applications, {\bf 424} (2007), 510-�539.

\bi{Dunkl} C.F.~Dunkl, {\it Integral kernels with reflection group
invariance}. Canadian Journal of Mathematics, {\bf 43} (1991)
1213--1227.

%\bi{Garna} A. El Garna, {\it The left-definite spectral theory for the
%Dunkl-Hermite differential-difference equation}, J. Math. Anal. Appl. {\bf 298} (2004) 463--486.

%\bi{Ev_Lit} W.N. Everitt, K.H. Kwon, L.L. Littlejohn, R.Wellman, {\it Orthogonal polynomial solutions of linear ordinary
%differential equations}, J. Comput. Appl. Math. {\bf 133} (2001), 85--109.

%\bi{Ger} Ya.L. Geronimus,\quad {\it Polynomials Orthogonal on a Circle
%and their Applications}, \\ Am.Math.Transl.,Ser.1 {\bf 3}(1962), 1-78.

%\bi{Ger1} Ya.L.Geronimus, {\it On polynomials orthogonal with respect to to
%the given numerical sequence and on Hahn's theorem}, Izv.Akad.Nauk, {\bf 4}
%(1940), 215-228 (in Russian).

\bi{Ger2} J.~Geronimus, {\it The orthogonality of some systems of
polynomials}, Duke Math. J. {\bf 14} (1947), 503--510.

\bi{GP} M.~F.~Gorodnii and G.~B.~Podkolzin, {\it Irreducible
representations of a graded Lie algebra}, in: Spectral theory of
operators and infinite-dimensional analysis, Institute of
Mathematics, Academy of Sciences of Ukraine, Kiev, (1984), 66–-77
(in Russian).

%\bi{Hahn} W.Hahn, {\it \"Uber die Jacobischen Polynome und Zwei
%verwandte Polynomklassen}, Math.Z. {\bf 39} (1935), 634-638.

%\bi{HN} E.Hendriksen and O. Nj{\aa}stad, {\it Biorthogonal Laurent polynomials
%with biorthogonal derivatives}, Rocky Mount. J. Math. {\bf 21} (1991),
%391-317.

%\bi{HR} E.Hendriksen and H. van Rossum, {\it Orthogonal Laurent polynomials},
%Indag. Math. (Ser. A) {\bf 48} (1986), 17-36.

\bi{Hild} F.B.~Hildebrand, {\it Introduction to Numerical
Analysis}. New York: McGraw-Hill, 1956.

%\bi{IsMas} M.E.H. Ismail and D. Masson, {\it Generalized orthogonality and
%continued fractions}, \\J.Approx.Theory {\bf 83} (1996), 1-40.

%\bi{JT} W.B.Jones and W.J.Thron, {\it Survey of continued fraction methods of
%solving moment problems} in: analytic Theory of Continued Fractions, LNM 932,
%Springer, Berlin, Heidelberg, New York (1981).

\bi{Ismail} M.E.H.~Ismail, {\it Classical and Quantum orthogonal
polynomials in one variable}. Encyclopedia of Mathematics and its
Applications (No. 98), Cambridge, 2005.

\bi{Ismail_Str} M.E.H.~Ismail, {\it Structure relations for
orthogonal polynomials}, Pacific J.Math. {\bf 240} (2009),
309--319.

\bi{KMG} S.~Karlin and J.~McGregor, {\it The Classification of
Birth and Death Processes}, Trans. Amer. Math. Soc., {\bf 86}
(1957), 366--400.

\bi{KLS} R.~Koekoek, P.~Lesky, R.~Swarttouw, {\it Hypergeometric
Orthogonal Polynomials and Their Q-analogues}, Springer-Verlag,
2010.

\bi{K_Str} T.~Koornwinder, {\it The structure relation for
Askey�Wilson polynomials}, J.Comp.Appl.Math., {\bf 207} (2007),
214--226.

\bi{KB} T.~Koornwinder and F.~Bouzeffour, {\it Nonsymmetric
Askey-Wilson polynomials as vector-valued polynomials}, Applicable
Analysis, {\bf 90} (2010), 731--746.  arXiv:1006.1140v1

\bi{K_Racah} T.~Koornwinder, {\it On the limit from q-Racah
polynomials to big q-Jacobi polynomials}, ArXiv:1011.5585.

\bi{Leonard} D.~Leonard, {\it Orthogonal Polynomials, Duality and
Association Schemes}, SIAM J. Math. Anal. {\bf 13} (1982)
656--663.

\bi{LR} L.L. ~Littlejohn and D.~ Race,  {\it Symmetric and  symmetrisable  differential  expressions},  Proc.  London Math.  Soc.  {\bf 60}  (1990)
344--364.

\bi{MP} F.~Marcell\'an and J.Petronilho, {\it Eigenproblems for
Tridiagonal 2-Toeplitz Matrices and Quadratic Polynomial
Mappings},  Lin. Alg. Appl. {\bf 260} (1997)   169--208.

\bi{NSU} A.F.~Nikiforov, S.K.~Suslov, and V.B.~Uvarov, {\em
Classical Orthogonal Polynomials of a Discrete Variable},
Springer, Berlin, 1991.

\bi{OS} V.L.Ostrovskyi and S.D.Silvestrov, {\it Representations of the real forms of a graded analogue
of the Lie algebra $sl(2, C)$}, Ukr. Mat. Zhurn. {\bf 44} (1992), no. 11, 1518--1524. English transl: Ukrain.Math.J. {\bf 44} (1992), 1395--1401.

\bi{RW} V.Rittenberg and D.Wyler, {\it Generalized superalgebras}, Nucl. Phys. {\bf B139} (1978), 189--202.

%\bi{Ros} M. Rosenblum, {\it Generalized Hermite Polynomials and the Bose-like Oscillator Calculus}, in: Oper. Theory
%Adv. Appl., vol. {\bf 73}, Birkhauser, Basel, 1994, pp. 369-�396.

\bi{ZhS_D}  V.~Spiridonov and A.~Zhedanov, {\it Discrete Darboux
transformations, the discrete-time Toda lattice, and the
Askey-Wilson polynomials}, Methods and Appl. Anal. , {\bf 2}
(1995), 369-�398.

%\bi{SVZ} V.Spiridonov, L.Vinet and A.Zhedanov, {\it Spectral
%transformations, self-similar reductions and orthogonal polynomials}, J.Phys.
%A:  Math.  and Gen.  {\bf 30} (1997), 7621--7637.

%\bi{Sz} G. Szeg\H{o}, Orthogonal Polynomials, fourth edition,  AMS, 1975.

\bi{TVZ_J} S.~Tsujimoto, L.~Vinet and A.~Zhedanov, {\it Jordan
algebras and orthogonal polynomials}, in preparation.

\bi{Ter} P.~Terwilliger, {\it Two linear transformations each
tridiagonal with respect to an eigenbasis of the other}. Linear
Algebra Appl. {\bf 330} (2001) 149--203.

\bi{Ter2} P.~Terwilliger, {\it Two linear transformations each
tridiagonal with respect to an eigenbasis of the other; an
algebraic approach to the Askey scheme of orthogonal polynomials},
arXiv:math/0408390v3.

%\bi{VZ_Bochner} L.Vinet and A.Zhedanov, {\it Generalized Bochner
%theorem: characterization of the Askey-Wilson polynomials},
%J.Comp.Appl.Math., {\bf 211} (2008) 45 -- 56.

\bi{Vidunas} R.~Vid\=unas, {\it Askey-Wilson relations and Leonard
pairs},  Discrete Mathematics {\bf 308} (2008) 479-–495.
arXiv:math/0511509v2.

\bi{VZ_little} L.~Vinet and A.~Zhedanov, {\it A ``missing`` family
of classical orthogonal polynomials}, J. Phys. A: Math. Theor.
{\bf 44} (2011) 085201,    arXiv:1011.1669v2.

\bi{VZ_big} L.~Vinet and A.~Zhedanov, {\it A limit $q=-1$ for big
q-Jacobi polynomials}, Trans.Amer.Math.Soc., to appear,
arXiv:1011.1429v3.

\bi{VZ_Bochner} L.~Vinet and A.~Zhedanov, {\it A Bochner theorem
for Dunkl polynomials}, SIGMA, {\bf 7} (2011), 020;
arXiv:1011.1457v3.

%\bi{WW} E.T. Whittacker, G.N. Watson, {\em A Course of Modern
%Analysis}, Cambridge, 1927.

\bi{Zhe} A.~Zhedanov, {\it "Hidden symmetry" of Askey-Wilson
polynomials}, Teoret. Mat. Fiz. {\bf 89} (1991) 190--204. (English
transl.: Theoret. and Math. Phys. {\bf 89} (1991), 1146--1157).

\bibitem{ZhS} A.~Zhedanov, {\it Rational spectral transformations
and orthogonal polynomials}, J. Comput. Appl. Math. {\bf 85}, no. 1
(1997), 67--86.

%\bi{ZheL} A.Zhedanov, {The "classical" Laurent biorthogonal
%polynomials}, J. Comp. Appl.Math. {\bf 98} (1998), 121--147.

\eb

\end{document}